\documentclass[a4paper,12pt]{amsart}

\numberwithin{equation}{section}

\DeclareMathOperator{\CAT}{CAT}
\DeclareMathOperator{\rad}{rad}

\DeclareMathOperator{\diam}{diam}

\DeclareMathOperator{\grad}{grad}
\DeclareMathOperator{\id}{id}

\DeclareMathOperator{\dis}{dis}
\DeclareMathOperator{\diag}{diag}
\DeclareMathOperator{\tr}{tr}
\DeclareMathOperator{\Min}{Min}
\DeclareMathOperator{\Fix}{Fix}

\theoremstyle{plain}
  \newtheorem{thm}{Theorem}[section]
  \newtheorem{lem}[thm]{Lemma}  
  \newtheorem{prop}[thm]{Proposition}
  \newtheorem{cor}[thm]{Corollary}

\theoremstyle{definition}
  
  \newtheorem{exmp}[thm]{Example}

  \newtheorem{rem}[thm]{Remark}

\begin{document}

\title
[Fixed point sets of parabolic isometries of CAT(0)-spaces]
{Fixed point sets of parabolic isometries\\ of CAT(0)-spaces}

\author
[K. Fujiwara]{Koji Fujiwara}
\author
[K. Nagano]{Koichi Nagano}
\author
[T. Shioya]{Takashi Shioya}

\thanks{The authors are partially supported by
the Grand-in-Aid for scientific research of Japanese Ministry of 
Education, Culture, Sports, and Technology.
The second author is partially supported by
the 2004 JSPS Postdoctoral Fellowships for Research Abroad.}

\email
[Koji Fujiwara]{fujiwara@math.tohoku.ac.jp}
\email
[Koichi Nagano]{nagano@math.tohoku.ac.jp}
\email
[Takashi Shioya]{shioya@math.tohoku.ac.jp}

\address
[Koji Fujiwara, Koichi Nagano, and Takashi Shioya]
{\endgraf Mathematical Institute, Tohoku University
\endgraf
Aoba, Sendai, Miyagi, 980-8578, Japan}

\date{\today}


\begin{abstract}
  We study the fixed point set in the ideal boundary of a parabolic
  isometry of a proper $\CAT(0)$-space.  We show that the radius of
  the fixed point set is at most $\pi/2$, and study its centers.  As a
  consequence, we prove that the set of fixed points is contractible
  with respect to the Tits topology.
\end{abstract}

\maketitle

\tableofcontents

\section
{Introduction}

$\CAT(0)$-spaces are generalizations of Hadamard
manifolds to geodesic spaces.
The classification of isometries of the hyperbolic
plane applies to isometries of $\CAT(0)$-spaces.

One of the nice results concerning hyperbolic
isometriesis is the flat torus theorem (cf.~\cite{ballmann,bridhaef}).
In the study of isometries, not only on $\CAT(0)$-spaces,
but also on Hadamard manifolds,
hyperbolic isometries have been more central than parabolic ones.
That would be explained from the view point of isometric group actions.
If a group acts cocompactly and
properly on a proper
$\CAT(0)$-space by isometries, then it does not contain
parabolic isometries.
However, if one considers an isometric group action which
are not cocompact, one may have to deal with parabolic
isometries.

In this paper, our focus is on parabolic isometries on $\CAT(0)$-spaces.
We study the fixed point sets of parabolic isometries in the ideal boundary
and generalize what Schroeder did
for parabolic isometries on Hadamard manifolds in Appendix 3 in
\cite{bagromsc}. But it is not straightforward because, for example,
much
less is known and available on analysis for $\CAT(0)$-spaces than
Hadamard manifolds. Also we need to treat certain $\CAT(1)$-spaces,
while they are only spheres in the argument for Hadamard manifolds in
\cite{bagromsc}.
At the end, as an example of our theorems, we examine a symmetric space
in detail.

\subsection
{Main theorems and examples}

Let $X$ be a complete $\CAT(0)$-space and $X(\infty)$ the ideal
boundary of $X$ defined as the asymptotic classes of rays in $X$.  We
classify an isometry $f$ of $X$ as elliptic, hyperbolic (axial), or
parabolic.  $f$ is called \emph{elliptic} if it has a fixed point in
$X$, and \emph{hyperbolic} if there exists a geodesic line (axis)
$\gamma$ in $X$ such that $f$ acts on $\gamma$ by a non-trivial
translation.  If $f$ is neither elliptic nor hyperbolic, then it is
called \emph{parabolic}.  We recall that $f$ is parabolic if and only
if the displacement function $d_f(p) := d(p,f(p))$ of $f$ does not
attain its minimum in $X$.  Any isometry of $X$ also acts on
$X(\infty)$ via rays.  It is known that if $X$ is proper (i.e., any
closed bounded subset is compact), then any parabolic isometry of $X$
has at least one fixed point in $X(\infty)$
(cf.~\cite{ballmann,bridhaef}).  In the case of improper $X$, there is
an example of a parabolic isometry $f$ of a separable Hilbert space
$X$ of infinite dimension such that $f$ has no fixed point in
$X(\infty)$ (and in $X$) (cf.~\cite{bridhaef}).  We denote by
$X_f(\infty)$ the fixed point set of $f$ in $X(\infty)$.

The ideal boundary $X(\infty)$ has a natural metric, called the
\emph{Tits metric $Td$}.  The metric space $(X(\infty),Td)$, say the
\emph{Tits ideal boundary}, is a complete $\CAT(1)$-space.
We define $\rad A := \inf_{x \in A} \sup_{y \in A} d(x,y)$ for a metric
space $A$ with metric $d$, which is called the \emph{radius of $A$}.
This notion for $A \subset X(\infty)$ is always defined for the Tits
metric.  For $p \in X$, we denote by $\Sigma_pX$ the space of
directions at $p$.  As one of the main results of this paper, we state
the following:

\begin{thm}\label{thm:radfix}
  Let $X$ be a proper $\CAT(0)$-space such that $\Sigma_pX$ is compact
  for every $p \in X$, and let $f$ be a parabolic isometry of $X$.
  Then we have $\rad X_f(\infty) \le \pi/2$.  In particular,
  $X_f(\infty)$ is contractible.
\end{thm}

Schroeder has proved Theorem \ref{thm:radfix} for smooth Hadamard
manifolds in Appendix 3 in \cite{bagromsc}.  In Theorem
\ref{thm:radfix}, the upper bound $\pi/2$ of $\rad X_f(\infty)$ is
optimal even for Hadamard manifolds (cf.~Example \ref{exmp:product}).
We also have some examples with $0 < \rad X_f(\infty) < \pi/2$
(cf.~Examples \ref{exmp:irreducible} and \ref{exmp:tiny}).  Notice
that if $X$ is proper and geodesically complete, then $\Sigma_pX$ is
compact for any $p \in X$.

Recall that $X$ is \emph{visible} if and only if $Td(x,y) = \infty$
for any distinct $x, y \in X(\infty)$.  By Theorem \ref{thm:radfix},
we immediately obtain:

\begin{cor}\label{cor:visible}
  Under the same assumption as in Theorem \ref{thm:radfix}, if $X$ is
  visible, then $X_f(\infty)$ consists of a single point.
\end{cor}

Buyalo \cite{sebuyalo} has shown that if $X$ is a complete, not
necessarily proper, Gromov-hyperbolic $\CAT(0)$-space, then $\inf d_f
= 0$, and $X_f(\infty)$ consists of a single point.  Let $X$ be a
proper $\CAT(0)$-space.  If $X$ is Gromov-hyperbolic, then $X$ is
visible.  If $X$ admits a cocompact group action, then the converse is
true (cf.~\cite{bridhaef}).

Next, we study the centers of $X_f(\infty)$.  A \emph{center of a
  metric space $A$} is defined to be a point in $A$ where the function
$A \ni x \mapsto \sup_{y \in A}d(x,y) \in [0,\infty]$ attains the
infimum, $\rad A$.  We denote by $C(A)$ the set of all centers of $A$,
and define $C^2(A) := C(C(A))$.

\begin{thm}\label{thm:dcenter}
  Let $X$ be a proper $\CAT(0)$-space of finite covering dimension
  such that $\Sigma_pX$ is compact for every $p \in X$.  Let $f$ be a
  parabolic isometry of $X$.  Then $C^2(X_f(\infty))$ consists of a
  single point, which is fixed by any isometry of $X$ leaving
  $X_f(\infty)$ invariant.  In particular, the point is a fixed point
  of any isometry of $X$ commuting $f$.
\end{thm}

Theorem \ref{thm:dcenter} for Hadamard manifolds has been shown by
Eberlein \cite{eberlein} following Schroeder's works in Appendix 3 in
\cite{bagromsc}.

We give some examples.

\begin{exmp}\label{exmp:product}
  Let us denote the hyperbolic plane by $\mathbb{H}^2$.  We consider
  the product Riemannian manifold $X := \mathbb{R} \times
  \underbrace{\mathbb{H}^2 \times \cdots \times
    \mathbb{H}^2}_{\text{$m$ times}}$, $m \ge 1$.  For $m$ parabolic
  isometries $h_1,h_2,\dots,h_m$ of $\mathbb{H}^2$, we define the
  product map $f := (\id_{\mathbb{R}}, h_1, \dots, h_m)$, where
  $\id_{\mathbb{R}}$ is the identity map on $\mathbb{R}$.  $f$ is a
  parabolic isometry of $X$.  We denote by $\mathbb{S}^{m-1}(1)$ the
  standard unit $(m-1)$-sphere in the Euclidean $m$-space
  $\mathbb{E}^m$ and set
  \begin{equation}
    \triangle_1^{m-1} := 
    \{ \, (x_1, \dots, x_m) \in \mathbb{S}^{m-1}(1) 
    \subset \mathbb{E}^m \mid 
    x_i \ge 0 \ \text{for all} \ i \, \},
    \label{eqn:sphspx}
  \end{equation}
  which we call the \emph{standard spherical $(m-1)$-simplex}.  We see
  that $X_f(\infty)$ is isometric to the spherical suspension over
  $\triangle_1^{m-1}$, where we refer \cite{burburiv} for the definition
  of spherical suspension.  We have $\rad X_f(\infty) = \pi/2$.
  $C(X_f(\infty))$ is isometric to $\triangle_1^{m-1}$ and
  $C^2(X_f(\infty))$ consists of the barycenter of $\triangle_1^{m-1}$.
\end{exmp}

\begin{exmp}\label{exmp:irreducible}
  We consider $X := SL(3,\mathbb{R})/SO(3)$, which is a
  five-dimensional, irreducible symmetric space of non-compact type
  and rank two.  $SL(3,\mathbb{R})$ is the identity component of the
  isometry group.  The Tits ideal boundary $(X(\infty),Td)$ is a thick
  spherical building of dimension one.  Weyl chambers of $X$ are
  corresponding to edges of the building $(X(\infty),Td)$ and any edge
  has length $\pi/3$.  According to our Theorem \ref{sl3}, for any
  parabolic isometry $f \in SL(3,\mathbb{R})$ of $X$, $X_f(\infty)$ is
  one of the following:
  \begin{enumerate}
  \item an edge,
  \item a closed interval of length $\pi$ consisting of three
    edges,
  \item the union of an edge $c$ and all edges incident to $c$.
  \end{enumerate}
  In (3), $X_f(\infty)$ has uncountably many edges.

  In Section 6, we precisely discuss all isometries in $SL(3,\mathbb{R})$.
\end{exmp}

For the irreducible symmetric space $X := SL(n,\mathbb{R})/SO(n)$, $n
\ge 3$, let $f$ be any isometry of $X$.  Since for any Weyl chamber
$c$ at infinity, $fc \cap c$ is a (possibly empty) face of $c$ and
since $\rad c \ge \pi/6$ (cf.~\cite{bridhaef}), we have either $\rad
X_f(\infty) = 0$ or $\ge \pi/6$.

For any given $\theta \in (0,\pi/2)$, we have an example with $\rad
X_f(\infty) = \theta$, where $X$ is a manifold with boundary.

\begin{exmp}\label{exmp:tiny}
  Let us take a parabolic isometry $h$ of $\mathbb{H}^2$ with its
  fixed point $x \in \mathbb{H}^2(\infty)$.  Let $\gamma$ be a ray in
  $\mathbb{H}^2$ tending to $x$, and $b_{\gamma}$ the Busemann
  function associated with $\gamma$ (see Section 2 for the definition
  of $b_\gamma$).  Note that $h$ leaves every horosphere
  $b_{\gamma}^{-1}(t)$ invariant.  For an arbitrarily given $\theta
  \in (0,\pi/2)$, we consider the closed convex subset
  \[
  X := \{ \; (p,s) \in \mathbb{H}^2 \times \mathbb{R} 
  \mid b_{\gamma}(p) \le -t, \, |s| \le t \sin \theta
  \ \text{for some $t \ge 0$}\;\}
  \]
  of $\mathbb{H}^2 \times \mathbb{R}$.  $X$ is a proper
  $\CAT(0)$-space and $(X(\infty),Td)$ is isometric to a closed
  interval of length $2\theta$ whose midpoint corresponds to $x$.  The
  product map $(h,\id_{\mathbb{R}})$ leaves $X$ invariant, and its
  restriction, say $f$, on $X$ is a parabolic isometry of $X$.
  Since $X_f(\infty)$ coincides with $X(\infty)$, we have $\rad
  X_f(\infty) = \theta$.
\end{exmp}

\subsection
{Key ideas of the proof of main theorems}

We prove Theorem \ref{thm:radfix} in Section 3 by using the
gradient curve for the displacement function, the existence of which
is established by Jost and Mayer \cite{jurgjost,uwemayer}.  The idea
of the proof of Theorem \ref{thm:radfix} is based on Schroeder's
original one for Hadamard manifolds in Appendix 3 in \cite{bagromsc}.
Since a $\CAT(0)$-space $X$ is not differentiable in general, we need
to investigate the directional derivatives of a Lipschitz continuous,
convex function on $X$.  It is non-trivial to prove a first variation
formula for such a function (see Lemma \ref{lem:1stvar}).

For Theorem \ref{thm:dcenter}, the original proof in \cite{bagromsc}
seems not to work for a $\CAT(0)$-space.  We find a new approach using
the geometry of the Tits ideal boundary $(X(\infty),Td)$ as explained
in the following.

For a topological space $Y$, we define a number $\dim_C Y$ as the
supremum of the covering dimensions of compact subsets of $Y$.  In
\cite{bkleiner} we see some geometric study for $\dim_C$ of a
$\CAT(0)$-space.  A key theorem to investigate the centers of
$X_f(\infty)$ is the following:

\begin{thm}\label{thm:smaller}
  Let $Y$ be a complete $\CAT(1)$-space of $\dim_C Y < \infty$ and
  diameter $\diam Y \le \pi/2$.  Then we have $\rad Y < \pi/2$.  In
  particular, $C(Y)$ consists of a single point.
\end{thm}

Schroeder has shown Theorem \ref{thm:smaller} if $Y$ is a closed
convex subset of a unit sphere in Appendix 3 in \cite{bagromsc}.  We
give a proof of Theorem \ref{thm:smaller} in Section 5.  The basic
strategy of the proof is following \cite{bagromsc}, however the
possible non-compactness of $Y$ makes the proof delicate much.  By the
non-compactness, we cannot avoid a discussion with error estimates,
which implies the stronger statement that there exists a constant
$\delta > 0$ depending only on $\dim_C Y < \infty$ such that $\rad Y
\le \pi/2 - \delta$.

It is necessary for Theorem \ref{thm:smaller} that $\dim_C Y$ is
finite. In fact, the inductive limit, say $Y$, of the standard
spherical $(m-1)$-simplices $\triangle_1^{m-1}$, $m = 1, 2, \dots$, in
\eqref{eqn:sphspx} is a complete $\CAT(1)$-space such that $\dim_C Y =
\infty$, $\diam Y = \pi/2$, and $\rad Y = \pi/2$.

For applying Theorem \ref{thm:smaller} to $Y := X_f(\infty)$, we need:

\begin{prop}\label{prop:cptdimension}
  For a proper $\CAT(0)$-space $X$ we have
  \[
  \dim_C (X(\infty),Td) \le \dim X - 1,
  \]
  where $\dim X$ denotes the covering dimension of $X$.
\end{prop}

Theorem C in \cite{bkleiner} implies Proposition
\ref{prop:cptdimension}, provided that $X$ has a cocompact group
action.  For the proof of the proposition, we rely on our result
in \cite{fujshiya} about the dimension of $X(\infty)$ with sphere topology.
There is an independent way to obtain the proposition using Lemma
11.1 of \cite{lytchak}.  We would like to thank A.~Lytchak for
bringing his work into our attention.
We do not know whether $\dim_C$ in
Proposition \ref{prop:cptdimension} can be replaced with the covering
dimension.

Theorem \ref{thm:dcenter} is proved in this way:  By Theorem
\ref{thm:radfix}, $Y := C(X_f(\infty))$ has $\diam Y \le \pi/2$.
Proposition \ref{prop:cptdimension} implies $\dim_C Y < \infty$.
Applying Theorem \ref{thm:smaller} leads to Theorem \ref{thm:dcenter}.
The details are stated in Section \ref{ssec:radfix}.

\section
{Preliminaries}

A \emph{minimizing geodesic} is, by definition, a length-minimizing
curve joining two points in a metric space.  We assume that all
minimizing geodesics have unit speed parameters.  Denote by
$\gamma_{pq}$ a minimizing geodesic from a point $p$ to a point $q$,
and by $[p,q]$ its image.  A \emph{geodesic triangle
  $\triangle(p,q,r)$} means a triple of minimizing geodesics
$\gamma_{pq}$, $\gamma_{qr}$, and $\gamma_{rp}$ for three points $p$,
$q$, and $r$, called \emph{vertices}.

For $\kappa \in \mathbb{R}$, let $M_{\kappa}^2$ be a complete, simply
connected model surface of constant curvature $\kappa$.  We set
$D_{\kappa} := \diam M_{\kappa}^2$.  Note that $D_{\kappa}$ is equal
to $\pi/\sqrt{\kappa}$ if $\kappa > 0$, and to $\infty$ if $\kappa \le
0$.  We say that a metric space $X$ is a \emph{$\CAT(\kappa)$-space}
if the following (1) and (2) are satisfied.
\begin{enumerate}
\item Any two points $p, q \in X$ with $d(p,q) < D_{\kappa}$ can be
  joined by a minimizing geodesic in $X$.
\item ($\CAT(\kappa)$-inequality) Let $\triangle(p,q,r)$ be any
  geodesic triangle in $X$ with perimeter $< 2D_{\kappa}$ and
  $\triangle(\widetilde{p},\widetilde{q},\widetilde{r})$ a comparison
  triangle of it in $M_\kappa^2$, i.e., having the same side lengths
  as $\triangle(\widetilde{p},\widetilde{q},\widetilde{r})$.  For any
  four points $x \in [p,q]$, $y \in [r,p]$, $\widetilde{x} \in
  [\widetilde{p},\widetilde{q}]$, and $\widetilde{y} \in
  [\widetilde{r},\widetilde{p}]$ such that $d(p,x) =
  d(\widetilde{p},\widetilde{x})$ and $d(p,y) =
  d(\widetilde{p},\widetilde{y})$, we have
  \[
  d(x,y) \le d(\widetilde{x},\widetilde{y}),
  \]
  where $d$ denotes the distance function.
\end{enumerate}

Let $X$ be a $\CAT(\kappa)$-space.  A minimizing geodesic
$\gamma_{pq}$ joining two points $p, q \in X$ with $d(p,q) <
D_{\kappa}$ is unique.  For $p \in X$ and $q_1, q_2 \in X \setminus
\{p\}$, we denote by $\angle_p(\gamma_{pq_1},\gamma_{pq_2})$ the angle
at $p$ between $\gamma_{pq_1}$ and $\gamma_{pq_2}$.  $\angle_p$ is a
pseudo-distance function on the set of all minimizing geodesics
emanating from $p$.  We denote by $\Sigma_p^*X$ its quotient
metric space by the relation $\angle_p=0$.  Let $\Sigma_pX$ be the
$\angle_p$-completion of $\Sigma_p^*X$, called the \emph{space of
  directions at $p$}.  Let $C_pX$ be the Euclidean cone over
$\Sigma_pX$, the \emph{tangent cone at $p$}.  $\Sigma_pX$ is a
complete $\mathrm{CAT}(1)$-space and $C_pX$ a complete
$\mathrm{CAT}(0)$-space.  We denote by $\dot{\gamma}(0)$ the
equivalence class in $\Sigma_p^*X$ represented by a minimizing
geodesic $\gamma$ from $p$.

Assume that $X$ is a complete $\CAT(0)$-space.  Two rays $\gamma,
\sigma : [0,+\infty) \to X$ are said to be \emph{asymptotic} if
$d(\gamma(t),\sigma(t))$ is uniformly bounded for all $t \ge 0$.  The
\emph{ideal boundary $X(\infty)$ of $X$} is defined as the set of all
asymptotic equivalence classes of rays in $X$.  $X(\infty)$ is
equipped with the sphere topology induced by the cone topology on $X
\sqcup X(\infty)$.  We denote by $\gamma(\infty)$ the equivalence
class in $X(\infty)$ represented by a ray $\gamma$ in $X$.  Notice
that for any $p \in X$ and $x \in X(\infty)$ there exists a unique ray
$\gamma_{px} \colon [0,\infty) \to X$ from $p$ to $\gamma(\infty) =
x$.  For $x, y \in X(\infty)$, we set $\angle(x,y) := \sup_{p \in X}
\angle_p(x,y)$, say the \emph{angle distance} between $x$ and $y$,
where we write $\angle_p(x,y) := \angle_p(\gamma_{px},\gamma_{py})$.
Note that $\angle$ is a distance function on $X(\infty)$ and is lower
semi-continuous with respect to the sphere topology.  We remark that
if $X$ is proper, then $X(\infty)$ is compact with respect to the
sphere topology.  The \emph{Tits distance on $X(\infty)$}, denoted by
$Td$, is the interior distance induced from $\angle$.  We have $\angle
= \min \{ Td, \pi \}$.  The Tits ideal boundary $(X(\infty),Td)$ of
$X$ is a complete $\CAT(1)$-space, which is non-compact in general.
The \emph{Busemann function $b_{\gamma} \colon X \to \mathbb{R}$
  associated with a ray $\gamma$ in $X$} is defined as
\[
b_{\gamma}(p) := \lim_{t \to \infty} \{d(p,\gamma(t)) - t\}.
\]
This is a $1$-Lipschitz continuous, convex function with
$b_{\gamma}(\gamma(0)) = 0$.

A subset $A$ of a metric space $X$ is said to be \emph{convex in $X$}
if any $x,y \in A$ can be joined by a minimizing geodesic and the
image of every such geodesic is contained in $A$.  If this condition
holds only for any $x, y \in A$ with $d(x,y) < r$, then $A$ is said to be
\emph{$r$-convex in $X$}.

Let $B$ be a closed subset of a metric space $X$.  We define a
function $d_B \colon X \to [0,\infty)$ by $d_B(p) := d(p,B)$, say the
\emph{distance function from $B$}.  For $p \in X \setminus B$, we
denote by $\gamma_{pB}$ a minimizing geodesic in $X$ from $p$ to $B$,
i.e., to a point $q \in B$ with $d_B(p) = d(p,q)$.

Assume that $B$ is a closed, convex subset of a complete
$\CAT(0)$-space.  Then, for any $p \in X$ there exists a unique point
$q \in B$ with $d_B(p) = d_B(p,q)$, in particular, $\gamma_{pq} =
\gamma_{pB}$.  We note that $d_B$ is a $1$-Lipschitz continuous,
convex function.

\section
{Estimate of radius of fixed point sets}

We prove Theorem \ref{thm:radfix}.

\subsection
{Directional derivatives of convex functions}

Let $X$ be a complete $\CAT(0)$-space and $F \colon X \to \mathbb{R}$
a locally Lipschitz continuous, convex function.  We discuss the
directional derivatives of $F$.  For any geodesic $\gamma$ in $X$, $F
\circ \gamma$ has left and the right derivatives.  Recall that the
tangent cone $C_pX$ is the quotient space $[0,+\infty) \times
\Sigma_pX/\{0\} \times \Sigma_pX$.  We identify the subspace $\{1\}
\times \Sigma_pX$ of $C_pX$ with $\Sigma_pX$.  Denote any element
$(t,v) \in C_pX$ by $tv$ and define $|tv| := t$.  Let $C_p^*X :=
[0,\infty) \times \Sigma_p^*X/\{0\} \times \Sigma_p^*X
\subset C_pX$.  The \emph{directional derivative $D_pF \colon
  C_p^*X \to \mathbb{R}$ of $F$ at a point $p \in X$} is defined
as
\[
D_pF(tv) := \lim_{s \to 0+}
\frac{F(\gamma_v(s))-F(\gamma_v(0))}{s} t,
\]
where $\gamma_v$ is a minimizing geodesic from $p$ with $v =
\dot{\gamma}_v(0)$.  The existence of the limit above is guaranteed by
the convexity of $F$.  $D_pF(tv)$ is independent of the choice of
$\gamma_v$.  $D_pF$ extends to a unique Lipschitz continuous function
on $C_pX$, which is convex (cf. Lemma 2.4 in \cite{bkleiner}).
Moreover, it is linear along each ray from the vertex $o_p$ of $C_pX$.

Assume that $\Sigma_pX$ is compact for every $p \in X$.  We say that a
point $p \in X$ is a \emph{critical point of $F$} if $D_pF(u) \ge 0$
for every $u \in \Sigma_pX$.  Note that, by the convexity of $F$, a
point is critical for $F$ if and only if it is a minimizer of $F$.
For more general functions, such as $c$-convex functions
(cf.~\cite{burburiv}), this is not true and we still have some local
properties stated below, e.g.~Theorem \ref{thm:gradcurve} and Lemma
\ref{lem:1stvar}.  By the convexity of $D_pF$ and the compactness of
$\Sigma_pX$, for any non-critical point $p$ of $F$, there exists a
unique direction $u_p \in \Sigma_pX$ where $D_pF |_{\Sigma_pX}$
attains its minimum ($< 0$).  We call $u_p$ the \emph{gradient
  direction of $-F$ at $p$}.  Define the \emph{gradient vector
  $\grad_p(-F) \in C_pX$ of $-F$ at a point $p$} by
\[
\grad_p(-F) := |D_pF(u_p)| u_p \in C_pX 
\]
if $p$ is non-critical, and by $\grad_p(-F) := o_p$ (the vertex) if
$p$ is critical.  It follows that $|\grad_p(-F)| = -D_pF(u_p)$.

\subsection
{Jost-Mayer's gradient curves}

We state a little restricted version of some theorems in
\cite{jurgjost,uwemayer}.

\begin{thm}[\cite{jurgjost,uwemayer}] \label{thm:gradcurve}
  Let $X$ be a complete $\CAT(0)$-space such that $\Sigma_pX$ is compact
  for every $p \in X$, and let $F \colon X \to \mathbb{R}$ be a convex
  function.  Then, for every $p \in X$ there exists a Lipschitz
  continuous curve $c_p \colon [0,\infty) \to X$ from $p = c_p(0)$,
  called the \emph{gradient curve from $p$ for $-F$}, such that for any
  $t \ge 0$ we have
  \begin{multline}
    \lim_{s \to 0+} \frac{d(c_p(t+s),c_p(t))}{s}
    = \lim_{s \to 0+} \frac{-F \circ c_p(t+s) + F \circ c_p(t)}
    {d(c_p(t+s),c_p(t))} \\    
    = \limsup_{q \to c_p(t)} \frac{-F(q)+F(c_p(t))}{d(q,c_p(t))}
    = |\grad_{c_p(t)}(-F)|,
    \tag{1}
  \end{multline}
  \begin{equation}
    (F \circ c_p)'_+(t) = |\grad_{c_p(t)}(-F)|^2,
    \tag{2}
  \end{equation}
  where $(F\circ c_p)'_+(t)$ is the right derivative of $F \circ c_p$
  at $t$.  Moreover, for any $r \ge 0$, the gradient curve
  $c_{c_p(t)}$ from $c_p(t)$ for $-F$ satisfies
  \[
  c_{c_p(t)}(r) = c_p(t+r).
  \]
\end{thm}

Under the same assumption as in Theorem \ref{thm:gradcurve},
we have:

\begin{lem}\label{lem:coincide}
  For the gradient curve $c_p$ from $p$ of $-F$, the right tangent
  vector $(\dot{c}_p)_+(0) \in C_pX$ exists and coincides with
  $\grad_p(-F)$.
\end{lem}

\begin{proof}
  By taking a sequence $\{s_i\}$ with $s_i \to 0+$, we have a limit $v
  \in \Sigma_pX$ of the direction $\dot{\gamma}_{pc_p(s_i)}(0)$ as $i
  \to \infty$.  By Theorem \ref{thm:gradcurve}(1), $D_pF(v)$ must be
  equal to $D_pF(u_p) = - |\grad_p(-F)|$.  We see that $v = u_p$ by
  the uniqueness of the gradient direction $u_p$.
\end{proof}

\subsection
{First variation formula}

The following is well-known.

\begin{lem}\label{lem:1stvarb}
  Let $X$ be a complete $\CAT(0)$-space. 
  \begin{enumerate}
  \item Let $B$ be a closed, convex subset of $X$.  Then for any $p
    \in X \setminus B$ and $v \in \Sigma_pX$ we have
    \[
    D_pd_B(v) = -\cos \angle_p(\dot{\gamma}_{pB}(0),v).
    \]
  \item Let $\gamma$ be a ray in $X$.  Then
    for any $p \in X$ and $v \in \Sigma_pX$ we have
    \[
      D_pb_{\gamma}(v) = 
      - \cos \angle_p (\dot{\gamma}_{p\gamma(\infty)}(0),v).
    \]
  \end{enumerate}
\end{lem}

\begin{proof}
  (1) follows from a standard argument (cf.~Section 4.5 of
  \cite{burburiv}).
  
  We prove (2).  Set $B_t := b_{\gamma}^{-1}(-\infty,-t]$ for $t > 0$.
  $B_t$ is convex in $X$.  Let $p \in X$ be any point.  If $t > 0$ is
  large enough for $p$, then $p \in X \setminus B_t$ and $d_{B_t}(p) =
  b_{\gamma}(p) + t$ (cf.~Proposition II.8.22 in \cite{bridhaef}),
  which and (1) imply (2).
\end{proof}

Some variants of Lemma \ref{lem:1stvarb}(1) is seen in Section 4.5 of
\cite{burburiv}. Note that the $\CAT(0)$-condition for $X$ is not
essential for Lemma \ref{lem:1stvarb}.

To prove a first variation formula for convex functions, we need
a lemma.

\begin{lem}\label{lem:sector}
  Let $S$ be a sector in $\mathbb{E}^2$ bounded by two distinct rays
  from the origin $o$.  Let $F \colon S \to \mathbb{R}$ be a function
  that is linear along each ray from $o$.  If the directional
  derivative $D_uF \colon C_uS \to \mathbb{R}$ of $F$ at a point $u
  \in S \setminus \{o\}$ exists, then $D_uF$ is linear on $C_uS$.
\end{lem}

Lemma \ref{lem:sector} is shown by a standard argument.  We omit the proof.

We prove the following first variation formula.

\begin{lem}\label{lem:1stvar}
  Let $F \colon X \to \mathbb{R}$ be a locally Lipschitz continuous,
  convex function on a complete $\CAT(0)$-space $X$.  Let $p \in X$ be
  a non-critical point of $F$ such that $\Sigma_pX$ is compact.  Then
  for any $v \in \Sigma_pX$ we have
  \[
  D_pF(v) \ge - |\grad_p(-F)| \cos \angle_p(u_p,v),
  \]
  where $u_p \in \Sigma_pX$ is the gradient direction of $-F$ at $p$.
\end{lem}

\begin{proof}
  Let $v \in \Sigma_pX$ be a direction.  If $\angle_p(u_p,v) = 0$, the
  lemma is obvious.  In the case where $\angle_p(u_p,v) = \pi$, the
  minimizing geodesic $\gamma_{u_pv}$ in $C_pX$ from $u_p$ to $v$
  passes through the vertex $o_p$, so that the convexity of $D_pF$
  along $\gamma_{u_pv}$ implies the lemma.
  
  We assume that $0 < \angle_p(u_p,v) < \pi$.  Consider the second
  derivative $D_{u_p}D_pF \colon C_{u_p}C_pX \to \mathbb{R}$.  Let $S
  \subset C_pX$ be the $2$-dimensional flat sector generated by
  $\gamma_{u_pv}$. $S$ is convex in $C_pX$.  We set $\xi :=
  \dot{\gamma}_{u_pv}(0)$ and $\eta := \dot{\gamma}_{u_po_p}(0)$, both
  which belong to $\Sigma_{u_p}S$.  Note that $C_{u_p}S$ is a flat
  half plane in $C_{u_p}C_pX$.  Take the direction $\zeta \in
  \Sigma_{u_p}S$ perpendicular to $\eta$.  Setting $\theta := \angle_{u_p}(o_p,v)$, we see
  $\xi = (\cos \theta)\eta + (\sin \theta)\zeta$.  Since Lemma
  \ref{lem:sector} implies the linearity of $D_{u_p}D_pF$, we have
  \[
  D_{u_p}D_pF(\xi) = D_{u_p}D_pF(\eta) \cos \theta
  + D_{u_p}D_pF(\zeta) \sin \theta.
  \]
  The linearity of $D_pF$ along $\gamma_{u_po_p}$ shows that
  $D_{u_p}D_pF(\eta) = - D_pF(u_p) > 0$.  Since $u_p$ is the minimum
  point of $D_pF$ on $\Sigma_pX$, we have $D_{u_p}D_pF(\zeta) \ge 0$.
  Thus, by noting $0 < \theta < \pi/2$,
  \begin{equation}
    \label{eq:1stvar1}
    D_{u_p}D_pF(\xi) \ge -D_pF(u_p) \cos\theta \ (> 0).
  \end{equation}
  It follows that the distance between $u_p$ and $v$ in $C_pX$ is
  equal to $2\cos\theta$, so that, by the convexity of $D_pF$ along
  $\gamma_{u_pv}$,
  \begin{equation}
    \label{eq:1stvar2}
    D_pF(v) \ge D_pF(u_p) + 2 \, D_{u_p}D_pF(\xi) \cos \theta.
  \end{equation}
  Combining \eqref{eq:1stvar1} and \eqref{eq:1stvar2} yields
  \[
  D_pF(v) \ge - D_pF(u_p) \cos 2\theta
  = D_pF(u_p) \cos \angle_p(u_p,v),
  \]
  which completes the  proof of Lemma \ref{lem:1stvar}.
\end{proof}

Note that the equality in Lemma \ref{lem:1stvar} does not
necessarily hold.  Lemma \ref{lem:1stvar} remains true for a locally
Lipschitz continuous, $c$-convex function $F$ on a locally
$\CAT(\kappa)$-space $X$, $c,\kappa \in \mathbb{R}$.

\subsection
{Monotone points} \label{ssec:mono}

Let $X$ be a complete $\CAT(0)$-space and $F \colon X \to \mathbb{R}$
a convex function.  The following terminology was introduced by
Eberlein in Section 4.1 of \cite{eberlein} for a Riemannian manifold.
A point $x \in X(\infty)$ is said to be \emph{$F$-monotone} if there
exists a ray $\gamma \colon [0,\infty) \to X$ with $x =
\gamma(\infty)$ such that $F \circ \gamma(t)$ is monotone
non-increasing in $t \ge 0$.  We denote by $X_F(\infty)$ the set of
all $F$-monotone points in $X(\infty)$, say the \emph{$F$-monotone
  set}.  For an isometry $f$ of $X$, we recall the displacement
function $d_f(p) := d(p,f(p))$, which is a $1$-Lipschitz continuous,
convex function on $X$.  For a ray $\gamma$ in $X$, $\gamma(\infty)$
is $d_f$-monotone if and only if $f\circ\gamma$ is asymptotic to
$\gamma$.  This leads to $X_{d_f}(\infty) = X_f(\infty)$.

The following lemma is obtained by the same discussion as in Section
4.1 of \cite{eberlein}.  We omit the proof.

\begin{lem}\label{lem:monoset}
  Let $F \colon X \to \mathbb{R}$ be a convex function.  Then we have
  the following (1), (2), and (3).
  \begin{enumerate}
  \item For a point $x \in X(\infty)$, the following (a),(b), and (c)
    are equivalent to each other.
    \begin{enumerate}
    \item[(a)] $x$ is $F$-monotone.
    \item[(b)] For any ray $\gamma$ with $x = \gamma(\infty)$, $F
      \circ \gamma(t)$ is monotone non-increasing in $t \ge 0$.
    \item[(c)] There exists a sequence $\{p_i\}$ of points in $X$
      converging to $x$ in the cone topology such that $F(p_i)$ is
      uniformly bounded from above.
    \end{enumerate}
  \item $X_F(\infty)$ is closed with respect to the sphere topology.
  \item If $X$ is proper, then $X_F(\infty)$ is a closed, $\pi$-convex
    subset of $(X(\infty),Td)$.
  \end{enumerate}
\end{lem}

\subsection
{Proof of Theorem \ref{thm:radfix}} \label{ssec:radfix} We prove the
theorem in the same way as in \cite{eberlein} by using Lemma
\ref{lem:1stvar}.  Let $X$ be a proper $\CAT(0)$-space such that
$\Sigma_pX$ is compact for every $p \in X$, and $f$ a parabolic
isometry of $X$.  Since the displacement function $d_f$ has no minimal
(or critical) point in $X$, we have the gradient direction $u_p$ of
$-d_f$ at any $p \in X$, which satisfies $D_pd_f(u_p) < 0$.  We fix a
point $p \in X$ and take the gradient curve $c_p$ from $p$ for $-d_f$.
By Lemma \ref{lem:coincide}, the right tangent vector
$(\dot{c}_p)_+(t) \in C_pX$ satisfies $(\dot{c}_p)_+(t) =
\grad_{c_p(t)}(-d_f)$ for any $t \ge 0$.  It follows from Theorem
\ref{thm:gradcurve}(1) that $d_f \circ c_p(t)$ is strictly monotone
decreasing in $t \ge 0$.  There exists a sequence $t_i \to \infty$
such that $c_p(t_i)$ converges to some point $x \in X(\infty)$ in the
cone topology.  Lemma \ref{lem:monoset}(1) implies $x \in X_f(\infty)$.
  
We take any $y \in X_f(\infty)$ and fix it.  It suffices to prove that
$Td(x,y) \le \pi/2$.  Let $v_t := \dot{\gamma}_{c_p(t)y}(0)$.
Consider the Busemann function $b := b_{\gamma_{py}}$ associated with
$\gamma_{py}$.  Since $y$ is $d_f$-monotone and by Theorem
\ref{thm:gradcurve}(1), Lemma \ref{lem:1stvarb}(2), and Lemma
\ref{lem:1stvar}, we have
\begin{align*}
  (b \circ c_p)_+'(t) &= - |\grad_{c_p(t)}(-d_f)| 
  \cos \angle_{c_p(t)}(u_{c_p(t)},v_t)\\
  &\le D_{c_p(t)}d_f(v_t) \le 0
\end{align*}
for any $t \ge 0$, and therefore $b \circ c_p(t)$ is monotone
non-increasing in $t$.  By Lemma \ref{lem:monoset}(1), $x$ is
$b$-monotone and, for any $q \in X$, $b \circ \gamma_{qx}(t)$ is
monotone non-increasing in $t$.  It follows from Lemma
\ref{lem:1stvarb}(2) that
\[
-\cos \angle_q (x,y) = 
(b \circ \gamma_{qx})'_+(0) \le 0,
\]
which proves $Td(x,y) \le \pi/2$.

Since $(X(\infty),Td)$ is $\CAT(1)$, $X_f(\infty)$ is contractible.
This completes the proof of Theorem \ref{thm:radfix}.  \qed

\medskip

Let $X$ be as in Theorem \ref{thm:radfix}.  Then we have $\rad
X_F(\infty) \le \pi/2$ for any locally Lipschitz continuous, convex
function $F$ on $X$ with no minimum in $X$.

\section
{Dimension of Tits ideal boundaries}

We need the following to prove Proposition \ref{prop:cptdimension}.

\begin{prop}[\cite{fujshiya}]\label{prop:cocov}
  Let $X$ be a proper $\CAT(0)$-space.  Then, the covering dimension
  of $X(\infty)$ for the sphere topology satisfies
  \[
  \dim X(\infty) \le \dim X - 1.
  \]
\end{prop}

\begin{proof}[Proof of Proposition \ref{prop:cptdimension}]
  By Proposition \ref{prop:cocov}, it suffices to show
  \[
  \dim_C (X(\infty),Td) \le \dim X(\infty).
  \]
  We consider the identity map $\iota \colon (X(\infty),Td) \to
  X(\infty)$, which is continuous.  Take any compact subset $K \subset
  (X(\infty),Td)$.  Since $X(\infty)$ is Hausdorff, $\iota|_K \colon K
  \to \iota(K)$ is a homeomorphism.  Thus, we have $\dim K = \dim
  \iota(K) \le \dim X(\infty)$.  This completes the proof.
\end{proof}

We denote by $\triangle^n = \triangle^n(a_0,a_1,\dots,a_n)$ a (closed)
$n$-simplex with vertices $a_0,a_1,\dots,a_n$.  Let $F_i \subset
\partial \triangle^n$ be the $(n-1)$-simplex that is the opposite face to
$a_i$, where $\partial \triangle^n$ is the boundary of $\triangle^n$. We say
that a map $\psi$ from $\triangle^n$ to a set \emph{collapses}
$\partial \triangle^n$ if
$$
\psi(F_0) \cap \psi(F_1) \cap \cdots \cap \psi(F_n) \neq \emptyset.
$$

The following is a consequence of Sperner's lemma (cf.~2.1 in
\cite{fedorchu}).

\begin{lem}\label{lem:sperner}
  Let $Y$ be a Hausdorff space of $\dim Y \le n-1$, $n \ge 1$.
  Then any continuous map $\psi : \triangle^n \to Y$
  collapses $\partial \triangle^n$.
\end{lem}

\begin{proof}
  Suppose that there exists a continuous map $\psi : \triangle^n \to Y$
  that does not collapse $\partial\triangle^n$.  We set $U_i := Y
  \setminus \psi(F_i)$, $i = 0, 1, \dots, n$, which are open in the
  Hausdorff space $Y$.  Since $\psi$ does not collapse $\partial
  \triangle^n$, $\{U_i\}_{i=0}^n$ is an open covering of $Y$.  By $\dim Y
  \le n-1$, there exists a refinement $\{V_i\}$ of $\{U_i\}$ of order
  at most $n$.  Since $\psi$ is continuous and the order of $\{V_i\}$
  is at most $n$, we can take a sufficiently refined triangulation of
  $\triangle^n$ such that for each simplex $s$ of it, $\psi(s)$
  intersects at most $n$ members of $\{V_i\}$.  Then we give a label
  by $i = 0, 1, \dots, n$ to each vertex of the refinement as follows.
  A label of a vertex $a$ is $i$ if $\psi(a) \in V_i$, which implies
  that this label is a Sperner label on $\triangle^n$.  Namely, each
  original vertex $a_i$ has the label $i$, and each vertex in the
  refinement contained in a $j$-dimensional simplex $\triangle^j =
  \triangle^j(a_{i_0},a_{i_1},\dots,a_{i_j})$ is labelled by one of
  $i_0,i_1,\dots,i_j$; e.g., a vertex contained in $F_i$ does not have
  the label $i$.  Therefore, by Sperner's lemma there exists at least
  one $n$-simplex $s^n$ in the refined triangulation of $\triangle^n$
  such that the vertices of $s^n$ have the $n+1$ different labels, $0,
  1, \dots, n$.
  
  On the other hand, since $\psi(s^n)$ is contained in at most $n$
  different $V_i$'s, the simplex $s^n$ has at most $n$ different
  labels.  This is a contradiction.
\end{proof}

Lemma \ref{lem:sperner} plays a key role in the proof of Theorem
\ref{thm:smaller} in Section 5.  As another application of Lemma
\ref{lem:sperner}, we have:

\begin{prop}\label{prop:geombdry}
  Let $Y$ be a $\CAT(1)$-space of $\dim_C Y \le m$, $m \ge 1$.  Then,
  for any embedding $\psi$ from an $m$-sphere $\mathbb{S}^m$ into $Y$
  we have $\rad \psi(\mathbb{S}^m) \ge \pi$.  In particular, if $m=1$,
  then $Y$ is locally an $\mathbb{R}$-tree.
\end{prop}

\begin{proof}
  Suppose that there exists an embedding $\psi : \mathbb{S}^m \to Y$
  satisfying $\rad \psi(\mathbb{S}^m) < \pi$.  Since $Y$ is $\CAT(1)$,
  $\psi(\mathbb{S}^m)$ is contractible in $Y$.  Hence, for a closed
  $(m+1)$-disk $D^{m+1}$ there is a continuous extension
  $\overline{\psi} \colon D^{m+1} \to Y$ of $\psi$.  By identifying
  $D^{m+1}$ with an $(m+1)$-simplex, $\overline{\psi}$ does not
  collapse $\partial D^{m+1}$ and $\dim \overline{\psi}(D^{m+1}) \le
  \dim_C Y \le m$.  This contradicts Lemma \ref{lem:sperner}.
\end{proof}

\begin{rem}\label{rem:remgeombdry}
  Let $X$ be a proper $\CAT(0)$-space of $\dim X \le n$.  By
  Proposition \ref{prop:cptdimension}, we can apply Proposition
  \ref{prop:geombdry} to $Y = (X(\infty),Td)$ and $m = n-1$.
\end{rem}

\section
{CAT(1)-spaces of small diameter}

We prove Theorems \ref{thm:dcenter} and \ref{thm:smaller}.

\subsection
{Small triangles}

Let $Y$ be a $\CAT(1)$-space.  For $x,y,z \in Y$ we set $\angle_x(y,z)
:= \angle_x(\gamma_{xy},\gamma_{xz})$.  Denote the image of
$\gamma_{xy}$ by $[x,y]$.  Let $\triangle = \triangle(a_0,a_1,a_2)$ be
a geodesic triangle in $Y$ with sides $[a_0,a_1],[a_1,a_2],[a_2,a_0]$,
and $\widetilde{\triangle} =
\triangle(\widetilde{a}_0,\widetilde{a}_1,\widetilde{a}_2)$ a
comparison triangle in $\mathbb{S}^2(1)$ of $\triangle$ with the same
side-lengths as $\triangle$.  Recall that $\angle_{a_i}(a_j,a_k) \le
\angle_{\widetilde{a}_i}(\widetilde{a}_j,\widetilde{a}_k)$ for
distinct $i, j, k = 0, 1, 2$.  We say that $\triangle(a_0,a_1,a_2)$ is
\emph{small} if $d(a_i,a_j) \le \pi/2$ for any $i, j = 0, 1, 2$.  If
$\triangle(a_0,a_1,a_2)$ is small, then we have $d(a_2,x) \le \pi/2$
for any $x \in [a_0,a_1]$ by the $\CAT(0)$-inequality.  If
$\triangle(a_0,a_1,a_2)$ is small and if $d(a_2,x) = \pi/2$ for some
$x \in [a_0,a_1] \setminus \{a_0,a_1\}$, then the triangle is an
isosceles triangle and bounds a convex spherical surface.

$O(\epsilon)$ denotes Landau's symbol, i.e., some universal function
such that $\limsup_{\epsilon \to 0} |O(\epsilon)|/\epsilon$ is finite.
We assume that $O(\epsilon)$ is positive.

For the proof of Theorem \ref{thm:smaller}, we first show:

\begin{lem}\label{lem:comp}
  Let $\epsilon \in (0,1)$ be a positive number.  Let $\triangle =
  \triangle(a_0,a_1,a_2)$ and $\triangle' = \triangle(a_0',a_1',a_2')$
  be small geodesic triangles in $Y$ and in $\mathbb{S}^2(1)$,
  respectively.  Then we have the following:
  \begin{enumerate}
  \item if $|d(a_i,a_j) - d(a_i',a_j')| \le \epsilon$ for any $i, j =
    0, 1, 2$ and if $d(a_0,a_j) \ge \epsilon^{1/2}$ for each $j = 1,
    2$, then we have
    \[
    \angle_{a_0}(a_1,a_2) < \angle_{a_0'}(a_1',a_2') + O(\epsilon^{1/2});
    \]
  \item if $\angle_{a_0}(a_1,a_2) \ge \angle_{a_0'}(a_1',a_2') -
    \epsilon$ and $|d(a_0,a_j) - d(a_0',a_j')| \le \epsilon$ for each
    $j = 1, 2$, then we have
    \[
    d(a_1,a_2) > d(a_1',a_2') - O(\epsilon).
    \]
  \end{enumerate}
\end{lem}

\begin{proof}
  (1): Let $\widetilde{\triangle} =
  \triangle(\widetilde{a}_0,\widetilde{a}_1,\widetilde{a}_2)$ be a
  comparison triangle in $\mathbb{S}^2(1)$ of $\triangle$.  Since $Y$
  is $\CAT(1)$, we have $\angle_{a_0}(a_1,a_2) \le
  \angle_{\widetilde{a}_0}(\widetilde{a}_1,\widetilde{a}_2)$.  By the
  assumption of $\triangle$ and $\triangle'$, we have
  the conclusion of (1).

  We omit the proof of (2).
\end{proof}

We next prove the following:

\begin{lem}\label{lem:fulltri}
  Let $\epsilon \in (0,1)$, and let $\triangle =
  \triangle(a_0,a_1,a_2)$ be a small geodesic triangle in $Y$.  Assume
  that there exists a point $y \in [a_0,a_1]$ such that
  $\min_{i=0,1}d(a_i,y) \ge \epsilon^{1/2}$ and
  $d(a_2,y) \ge \pi/2 - \epsilon$.  Then we have
  \begin{gather}
    |\;\angle_y(a_2,a_i) - \pi/2\;| < O(\epsilon^{1/2}),
    \quad i = 0,1,\tag{1}\\
    d(a_2,x) > \pi/2 - O(\epsilon^{1/2})\tag{2}
  \end{gather}
  for any $x \in [a_0,a_1]$.
\end{lem}

\begin{proof}
  (1): Let $\triangle_i' = \triangle(y',a_i',a_2')$, $i = 0, 1$, be
  two spherical triangles in $\mathbb{S}^2(1)$ such that $d(y',a_i') =
  d(y,a_i)$, $d(a_i',a_2') = d(a_i,a_2)$, and $d(a_2',y') = \pi/2$.
  Since each $\triangle_i'$ is small, we have $\angle_{y'}(a_2',a_i')
  \le \pi/2$.  By $d(a_2,y) \ge \pi/2 - \epsilon$, we have $|d(a_2,y)
  - d(a_2',y')| \le \epsilon$.  Applying Lemma \ref{lem:comp}(1) to
  $\triangle(y,a_i,a_2)$ and $\triangle_i'$ yields that
  $\angle_y(a_2,a_i) < \pi/2 + O(\epsilon^{1/2})$.  Therefore, by $\pi
  \le \angle_y(a_2,a_0) + \angle_y(a_2,a_1)$ we have
  $\angle_y(a_2,a_i) > \pi/2 - O(\epsilon^{1/2})$.
  
  (2): For any given $x \in [a_0,a_1] \setminus \{y\}$, let us take a
  small spherical isosceles triangle $\triangle'' =
  \triangle(y'',x'',a_2'')$ such that $d(a_2'',x'') = d(a_2'',y'') =
  \pi/2$ and $d(x'',y'') = d(x,y)$.  Since
  $\angle_{y''}(a_2'',x'') = \pi/2$ and by (1) we have
  $\angle_y(a_2,x) > \angle_{y''}(a_2'',x'') - O(\epsilon^{1/2})$.
  Applying Lemma \ref{lem:comp}(2) to $\triangle(y,x,a_2)$ and
  $\triangle''$ shows (2).
\end{proof}

\begin{lem}\label{lem:ruled}
  Let $\epsilon \in (0,1)$ be a positive number.  For two small
  geodesic triangles $\triangle = \triangle(a_0,a_1,a_2)$ in $Y$ and
  $\triangle' = \triangle(a_0',a_1',a_2')$ in $\mathbb{S}^2(1)$, we
  assume that
  \begin{enumerate}
  \item $d(a_2,x) > \pi/2 - \epsilon$ for any $x
    \in [a_0,a_1]$;
  \item $d(a_2',x') = \pi/2$ for any $x' \in
    [a_0',a_1']$;
  \item $|d(a_0,a_1) - d(a_0',a_1')| < \epsilon$.
  \end{enumerate}
  For any four points $x_i \in [a_2,a_i]$, $x_i' \in [a_2',a_i']$, $i
  = 0, 1$, such that $d(a_2,x_i)/d(a_2,a_i) =
  d(a_2',x_i')/d(a_2',a_i')$, we have
  \begin{equation}
    |d(x_0,x_1) - d(x_0',x_1')| < O(\epsilon^{1/4}).
  \end{equation}
\end{lem}

\begin{proof}
  Take such four points $x_0$, $x_1$, $x_0'$, and $x_1'$.
  We may assume that $d(a_0',a_1') \ge 4\epsilon^{1/2}$.  Note that
  $d(a_0,a_1) \ge 3\epsilon^{1/2}$.  Take $\overline{y}_i \in
  [a_0,a_1]$ and $\overline{y}_i' \in [a_0',a_1']$ with
  $d(a_i,\overline{y}_i) = d(a_i',\overline{y}_i') = \epsilon^{1/2}$
  for $i = 0, 1$.  Let $y_i \in [a_2,\overline{y}_i]$ and $y_i' \in
  [a_2',\overline{y}_i']$ be the points determined by
  $$
  \frac{d(a_2,y_i)}{d(a_2,\overline{y}_i)} = \frac{d(a_2,x_i)}{d(a_2,a_i)},
  \quad
  \frac{d(a_2',y_i')}{d(a_2',\overline{y}_i')}
  = \frac{d(a_2',x_i')}{d(a_2',a_i')}
  $$
  (cf.~Figure \ref{fig:ruledp}).  Let $\widetilde{\triangle} =
  \triangle(\widetilde{\overline{y}}_0,\widetilde{\overline{y}}_1,\widetilde{a}_2)$
  be a spherical comparison triangle in $\mathbb{S}^2(1)$ of
  $\triangle = \triangle(\overline{y}_0,\overline{y}_1,a_2)$, and
  $\widetilde{y}_0, \widetilde{y}_1 \in \widetilde{\triangle}$ the
  corresponding points to $y_0, y_1$.  Considering the two geodesic
  triangles $\widetilde{\triangle}$ and
  $\triangle(\overline{y}_0',\overline{y}_1',a_2')$ in
  $\mathbb{S}^2(1)$, we have
  $|d(\widetilde{y}_0,\widetilde{\overline{y}}_1) -
  d(y_0',\overline{y}_1')|, |d(\widetilde{y}_0,\widetilde{y}_1) -
  d(y_0',y_1')| < O(\epsilon)$,
  \begin{equation}
    d(y_0,\overline{y}_1) < d(y_0',\overline{y}_1') + O(\epsilon),
    \label{eqn:ruled0}
  \end{equation}
  and $d(y_0,y_1) < d(y_0',y_1') + O(\epsilon)$.  By $d(x_i,y_i),
  d(x_i',y_i') < O(\epsilon^{1/2})$, we have $d(x_0,x_1) <
  d(x_0',x_1') + O(\epsilon^{1/2})$.  To obtain the opposite
  inequality, it suffices to prove
  \begin{equation}
    d(y_0,y_1) > d(y_0',y_1') - O(\epsilon^{1/4}).
    \label{eqn:ruled1}
  \end{equation}

\begin{figure}[t]
\unitlength 0.1in
\begin{picture}(48.50,19.30)(6.30,-22.30)
%
\special{pn 8}%
\special{pa 800 2200}%
\special{pa 2790 2200}%
\special{fp}%
%
\special{pn 8}%
\special{pa 1800 610}%
\special{pa 810 2200}%
\special{fp}%
%
\special{pn 8}%
\special{pa 1800 600}%
\special{pa 2800 2200}%
\special{fp}%
%
\special{pn 8}%
\special{pa 1800 620}%
\special{pa 1200 2200}%
\special{fp}%
%
\special{pn 8}%
\special{pa 1800 620}%
\special{pa 2390 2200}%
\special{fp}%
%
\special{pn 8}%
\special{sh 1}%
\special{ar 1800 610 10 10 0  6.28318530717959E+0000}%
\special{sh 1}%
\special{ar 1780 610 10 10 0  6.28318530717959E+0000}%
%
\special{pn 8}%
\special{sh 1}%
\special{ar 810 2200 10 10 0  6.28318530717959E+0000}%
\special{sh 1}%
\special{ar 810 2200 10 10 0  6.28318530717959E+0000}%
%
\special{pn 8}%
\special{sh 1}%
\special{ar 2800 2210 10 10 0  6.28318530717959E+0000}%
\special{sh 1}%
\special{ar 2800 2210 10 10 0  6.28318530717959E+0000}%
%
\special{pn 8}%
\special{sh 1}%
\special{ar 1180 2200 10 10 0  6.28318530717959E+0000}%
\special{sh 1}%
\special{ar 1180 2200 10 10 0  6.28318530717959E+0000}%
%
\special{pn 8}%
\special{sh 1}%
\special{ar 2390 2190 10 10 0  6.28318530717959E+0000}%
\special{sh 1}%
\special{ar 2390 2190 10 10 0  6.28318530717959E+0000}%
%
\special{pn 8}%
\special{sh 1}%
\special{ar 1190 1600 10 10 0  6.28318530717959E+0000}%
\special{sh 1}%
\special{ar 1190 1600 10 10 0  6.28318530717959E+0000}%
%
\special{pn 8}%
\special{sh 1}%
\special{ar 1390 1690 10 10 0  6.28318530717959E+0000}%
\special{sh 1}%
\special{ar 1390 1690 10 10 0  6.28318530717959E+0000}%
%
\special{pn 8}%
\special{sh 1}%
\special{ar 2100 1450 10 10 0  6.28318530717959E+0000}%
\special{sh 1}%
\special{ar 2100 1450 10 10 0  6.28318530717959E+0000}%
%
\special{pn 8}%
\special{sh 1}%
\special{ar 2280 1360 10 10 0  6.28318530717959E+0000}%
\special{sh 1}%
\special{ar 2280 1360 10 10 0  6.28318530717959E+0000}%
%
\special{pn 8}%
\special{pa 1400 1700}%
\special{pa 2390 2200}%
\special{fp}%
%
\special{pn 8}%
\special{pa 1390 1700}%
\special{pa 2100 1460}%
\special{fp}%
\put(17.5000,-4.7000){\makebox(0,0)[lb]{$a_2$}}%
\put(12.2000,-17.1000){\makebox(0,0)[lb]{$y_0$}}%
\put(22.2000,-15.6000){\makebox(0,0)[lb]{$y_1$}}%
\put(23.9000,-12.7000){\makebox(0,0)[lb]{$x_1$}}%
\put(9.9000,-14.8000){\makebox(0,0)[lb]{$x_0$}}%
\put(11.2000,-24.0000){\makebox(0,0)[lb]{$\overline{y}_0$}}%
\put(23.1000,-23.9000){\makebox(0,0)[lb]{$\overline{y}_1$}}%
\put(6.3000,-24.0000){\makebox(0,0)[lb]{$a_0$}}%
\put(29.0000,-24.0000){\makebox(0,0)[lb]{$a_1$}}%
%
\special{pn 8}%
\special{pa 3380 2200}%
\special{pa 5370 2200}%
\special{fp}%
%
\special{pn 8}%
\special{pa 4380 610}%
\special{pa 3390 2200}%
\special{fp}%
%
\special{pn 8}%
\special{pa 4380 600}%
\special{pa 5380 2200}%
\special{fp}%
%
\special{pn 8}%
\special{pa 4380 620}%
\special{pa 3780 2200}%
\special{fp}%
%
\special{pn 8}%
\special{pa 4380 620}%
\special{pa 4970 2200}%
\special{fp}%
%
\special{pn 8}%
\special{sh 1}%
\special{ar 4380 610 10 10 0  6.28318530717959E+0000}%
\special{sh 1}%
\special{ar 4360 610 10 10 0  6.28318530717959E+0000}%
%
\special{pn 8}%
\special{sh 1}%
\special{ar 3390 2200 10 10 0  6.28318530717959E+0000}%
\special{sh 1}%
\special{ar 3390 2200 10 10 0  6.28318530717959E+0000}%
%
\special{pn 8}%
\special{sh 1}%
\special{ar 5380 2210 10 10 0  6.28318530717959E+0000}%
\special{sh 1}%
\special{ar 5380 2210 10 10 0  6.28318530717959E+0000}%
%
\special{pn 8}%
\special{sh 1}%
\special{ar 3760 2200 10 10 0  6.28318530717959E+0000}%
\special{sh 1}%
\special{ar 3760 2200 10 10 0  6.28318530717959E+0000}%
%
\special{pn 8}%
\special{sh 1}%
\special{ar 4970 2190 10 10 0  6.28318530717959E+0000}%
\special{sh 1}%
\special{ar 4970 2190 10 10 0  6.28318530717959E+0000}%
%
\special{pn 8}%
\special{sh 1}%
\special{ar 3770 1600 10 10 0  6.28318530717959E+0000}%
\special{sh 1}%
\special{ar 3770 1600 10 10 0  6.28318530717959E+0000}%
%
\special{pn 8}%
\special{sh 1}%
\special{ar 3970 1690 10 10 0  6.28318530717959E+0000}%
\special{sh 1}%
\special{ar 3970 1690 10 10 0  6.28318530717959E+0000}%
%
\special{pn 8}%
\special{sh 1}%
\special{ar 4680 1450 10 10 0  6.28318530717959E+0000}%
\special{sh 1}%
\special{ar 4680 1450 10 10 0  6.28318530717959E+0000}%
%
\special{pn 8}%
\special{sh 1}%
\special{ar 4860 1360 10 10 0  6.28318530717959E+0000}%
\special{sh 1}%
\special{ar 4860 1360 10 10 0  6.28318530717959E+0000}%
%
\special{pn 8}%
\special{pa 3980 1700}%
\special{pa 4970 2200}%
\special{fp}%
%
\special{pn 8}%
\special{pa 3970 1700}%
\special{pa 4680 1460}%
\special{fp}%
\put(43.3000,-4.7000){\makebox(0,0)[lb]{$a_2'$}}%
\put(38.0000,-17.1000){\makebox(0,0)[lb]{$y_0'$}}%
\put(48.0000,-15.6000){\makebox(0,0)[lb]{$y_1'$}}%
\put(49.7000,-12.7000){\makebox(0,0)[lb]{$x_1'$}}%
\put(35.7000,-14.8000){\makebox(0,0)[lb]{$x_0'$}}%
\put(37.0000,-24.0000){\makebox(0,0)[lb]{$\overline{y}_0'$}}%
\put(48.9000,-23.9000){\makebox(0,0)[lb]{$\overline{y}_1'$}}%
\put(32.1000,-24.0000){\makebox(0,0)[lb]{$a_0'$}}%
\put(54.8000,-24.0000){\makebox(0,0)[lb]{$a_1'$}}%
\end{picture}%
\caption{$\triangle = \triangle(a_0,a_1,a_2)$ and 
$\triangle' = \triangle(a_0',a_1',a_2')$}
\label{fig:ruledp}
  \end{figure}
  
  Applying Lemma \ref{lem:fulltri}(1) to
  $\triangle(a_0,\overline{y}_1,a_2)$ and
  $\triangle(\overline{y}_0,a_1,a_2)$ yields
  \begin{equation}
    \pi/2 - O(\epsilon^{1/2}) <
    \angle_{\overline{y}_0}(a_2,\overline{y}_1), \,
    \angle_{\overline{y}_1}(a_2,\overline{y}_0)
    < \pi/2 + O(\epsilon^{1/2}).
    \label{eqn:ruleda}
  \end{equation}
  Consider $\triangle(\overline{y}_0,\overline{y}_1,y_0)$ and
  $\triangle(\overline{y}_0',\overline{y}_1',y_0')$.  By
  $\angle_{\overline{y}_0'}(a_2',\overline{y}_1') = \pi/2$ and
  \eqref{eqn:ruleda}, we have
  $\angle_{\overline{y}_0}(a_2,\overline{y}_1) >
  \angle_{\overline{y}_0'}(a_2',\overline{y}_1') - O(\epsilon^{1/2})$.
  Hence Lemma \ref{lem:comp}(2) implies $d(y_0,\overline{y}_1) >
  d(y_0',\overline{y}_1') - O(\epsilon^{1/2})$.  This together with
  \eqref{eqn:ruled0} implies
  \begin{equation}
    |d(y_0,\overline{y}_1) - d(y_0',\overline{y}_1')| < O(\epsilon^{1/2}).
    \label{eqn:ruledb}
  \end{equation}
  Therefore, by Lemma \ref{lem:comp}(1) we see that
  $\angle_{\overline{y}_1}(\overline{y}_0,y_0) <
  \angle_{\overline{y}_1'}(\overline{y}_0',y_0') + O(\epsilon^{1/4})$.
  It follows from \eqref{eqn:ruleda} and
  $\angle_{\overline{y}_1'}(a_2',\overline{y}_0') = \pi/2$ that
  \begin{multline}
    \angle_{\overline{y}_1}(y_0,y_1) 
    \ge \angle_{\overline{y}_1}(y_1,\overline{y}_0) - \angle_{\overline{y}_1}(\overline{y}_0,y_0) \\
    > \pi/2 - \angle_{\overline{y}_1'}(\overline{y}_0',y_0') - O(\epsilon^{1/4}) 
    = \angle_{\overline{y}_1'}(y_0',y_1') - O(\epsilon^{1/4}).
    \label{eqn:ruledc}
  \end{multline}
  By \eqref{eqn:ruledb}, \eqref{eqn:ruledc}, and Lemma
  \ref{lem:comp}(2), we have \eqref{eqn:ruled1}.  This completes the
  proof.
\end{proof}

\subsection
{Proof of Theorem \ref{thm:smaller}}

We need a lemma.

\begin{lem}\label{lem:discurve}
  Let $\epsilon$ and $l$ be positive numbers, and let $c : [0,l] \to Y$
  be a $1$-Lipschitz continuous curve from a point $x_0$ to a point
  $x_1$ in a metric space $Y$ such that
  \begin{equation}
    l  < d(x_0,x_1) + \epsilon.
    \label{eqn:discurve1}
  \end{equation}
  Assume that there exists a minimizing geodesic $\gamma_{x_0x_1}$
  joining $x_0$ to $x_1$.  Then, for any $s \in [0,1]$, setting $x_s
  := \gamma_{x_0x_1}(s\,d(x_0,x_1))$ we have
  \[
  d(x_s,c(sl)) < 2\epsilon.
  \]
\end{lem}

Note that the parameter of $c$ is not necessarily proportional to
the arc-length.

\begin{proof}
  Since $c$ is $1$-Lipschitz continuous, it follows from
  \eqref{eqn:discurve1} that
  $$
  d(x_0,c(sl)) + d(c(sl),x_1) \le sl + (1-s)l < d(x_0,x_1) + \epsilon,
  $$
  and hence, by the triangle inequality,
  \begin{equation}
    0 \le sl - d(x_0,c(sl)) < \epsilon.
    \label{eqn:discurvea}
  \end{equation}
  By \eqref{eqn:discurve1} and $d(x_s,x_1) \le (1-s)l$, we have
  $$
  sl \ge d(x_0,x_s) = d(x_0,x_1) - d(x_s,x_1) > sl - \epsilon.
  $$
  Combining this and \eqref{eqn:discurvea} yields
  \[
  |d(x_0,x_s) - d(x_0,c(sl))| < 2\epsilon.
  \]
  By the triangle inequality, this completes the proof.
\end{proof}

Let $Y$ be a $\CAT(1)$-space with $\diam Y \le \pi/2$, and let
$\rho \colon Y \to \mathbb{R}$ be the function defined by $\rho(x) :=
\sup_{y \in Y} d(x,y)$.  By the definition, $\rad Y = \inf_{x \in Y}
\rho(x) \le \pi/2$.  We define the constant $\delta_m := \pi/2 - \rad
\triangle_1^m$, where $\triangle_1^m$ is the standard spherical simplex
defined in \eqref{eqn:sphspx}.  $\delta_m$ is strictly monotone
decreasing in $m = 1, 2, \dots$.  Denote the barycenter of
$\triangle_1^m$ by $b_m'$.

The \emph{distortion $\dis\varphi$} of a map $\varphi \colon A_1 \to
A_2$ between metric spaces is defined by
$$
\dis \varphi := \sup_{x,y \in A_1} | d(\varphi(x),\varphi(y)) - d(x,y) |.
$$

We prove the following:

\begin{lem}\label{lem:induction}
  Let $\epsilon$ be a positive number with $\epsilon \ll \delta_m$.
  Assume that there exists a $1$-Lipschitz continuous map $\varphi_m
  \colon \triangle_1^m \to Y$ such that $\dis \varphi_m < \epsilon$
  and $\rho(b_m) > \pi/2-\epsilon$, where $b_m := \varphi_m(b_m')$.
  Then, there exists a $1$-Lipschitz continuous map $\varphi_{m+1}
  \colon \triangle_1^{m+1} \to Y$ such that $\dis \varphi_{m+1} <
  O(\epsilon^{1/8})$.
\end{lem}

\begin{proof}
  Denote by $a_0',\dots,a_{m+1}'$ the vertices of $\triangle_1^{m+1}$,
  and set $a_i := \varphi_m(a_i')$.  Let $\triangle_1^m \subset
  \partial \triangle_1^{m+1}$ be the face opposite to $a_{m+1}'$.
  There exists a point $a_{m+1} \in Y$ with $d(a_{m+1},b_m) > \pi/2 -
  \epsilon$.  We construct a map $\varphi_{m+1} \colon
  \triangle_1^{m+1} \to Y$ as follows.  For any given $x' \in
  \triangle_1^{m+1}$, the segment $[a_{m+1}',x']$ extends to a segment
  $[a_{m+1}',\overline{x}']$ with $\overline{x}' \in \triangle_1^m$.
  Set $\overline{x} := \varphi_m(\overline{x}')$.  There is a unique
  point $x \in [a_{m+1},\overline{x}]$ such that
  \[
  \frac{d(a_{m+1},x)}{d(a_{m+1},\overline{x})}
  = \frac{d(a_{m+1}',x')}{d(a_{m+1}',\overline{x}')}.
  \]
  We then define $\varphi_{m+1}(x') := x$.  It follows that
  $\varphi_{m+1}(a_{m+1}') = a_{m+1}$ and $\varphi_{m+1}|_{\triangle_1^m}
  = \varphi_m$.  Note that $\varphi_m$ and $\varphi_{m+1}$ are not
  necessarily injective.
  
  Let us prove that for any $z \in \varphi_m(\triangle_1^m)$,
  \begin{equation}
    d(a_{m+1},z) > \pi/2 - O(\epsilon^{1/2}).
    \label{eqn:smallera}
  \end{equation}
  Take a point $z' \in \triangle_1^m$ with $\varphi_m(z') = z$.  The
  segment $[b_m',z']$ extends to a segment $[z_0',z_1']$ with
  $z_0',z_1' \in \partial \triangle_1^m$.  Since $\delta_m$ coincides
  with the radius of the inscribed sphere of $\triangle_1^m$ centered at
  $b_m'$, we have $d(b_m',z_i') \ge \delta_m$ for each $i = 0, 1$.
  Set $z_i := \varphi_m(z_i')$.  Consider the $1$-Lipschitz continuous
  curve $c := \varphi_m \circ \gamma_{z_0'z_1'}$ joining $z_0$ and
  $z_1$.  Note that $c$ passes through $z$ and $b_m$.  Choose a number
  $s \in [0,1]$ with $c(s\,d(z_0',z_1')) = b_m$ and let $b :=
  \gamma_{z_0z_1}(s\,d(z_0,z_1))$ (cf.~Figure \ref{fig:smallerp}).
  Since $\dis \varphi_m < \epsilon$, we see that $d(z_0',z_1') <
  d(z_0,z_1) + \epsilon$.  Lemma \ref{lem:discurve} implies that
  $d(b,b_m) < 2\epsilon$ and so $d(a_{m+1},b) > \pi/2 - 3\epsilon$ by
  the assumption for $a_{m+1}$.  By $\epsilon \ll \delta_m$ we have
  $$
  d(b,z_i) > d(b_m',z_i') - 3\epsilon \ge \delta_m - 3\epsilon
  \ge \epsilon^{1/2}
  $$
  for each $i = 0,1$.
  Applying Lemma \ref{lem:fulltri}(2) to
  $\triangle(z_0,z_1,a_{m+1})$ yields that $d(a_{m+1},y) > \pi/2 -
  O(\epsilon^{1/2})$ for any $y \in [z_0,z_1]$.  Therefore, by Lemma
  \ref{lem:discurve} we obtain \eqref{eqn:smallera}.

\begin{figure}
\unitlength 0.1in
\begin{picture}(48.00,20.90)(3.10,-23.40)
%
\special{pn 8}%
\special{pa 610 2010}%
\special{pa 2400 2010}%
\special{fp}%
%
\special{pn 8}%
\special{pa 1490 590}%
\special{pa 610 2010}%
\special{fp}%
%
\special{pn 8}%
\special{pa 1500 590}%
\special{pa 2390 2020}%
\special{fp}%
%
\special{pn 8}%
\special{sh 1}%
\special{ar 1500 2010 10 10 0  6.28318530717959E+0000}%
\special{sh 1}%
\special{ar 1500 2010 10 10 0  6.28318530717959E+0000}%
%
\special{pn 8}%
\special{sh 1}%
\special{ar 1490 600 10 10 0  6.28318530717959E+0000}%
\special{sh 1}%
\special{ar 1490 600 10 10 0  6.28318530717959E+0000}%
%
\special{pn 8}%
\special{sh 1}%
\special{ar 630 2010 10 10 0  6.28318530717959E+0000}%
\special{sh 1}%
\special{ar 630 2010 10 10 0  6.28318530717959E+0000}%
%
\special{pn 8}%
\special{sh 1}%
\special{ar 2400 2010 10 10 0  6.28318530717959E+0000}%
\special{sh 1}%
\special{ar 2400 2010 10 10 0  6.28318530717959E+0000}%
%
\special{pn 8}%
\special{pa 3160 2010}%
\special{pa 4950 2010}%
\special{fp}%
%
\special{pn 8}%
\special{pa 4040 590}%
\special{pa 3160 2010}%
\special{fp}%
%
\special{pn 8}%
\special{pa 4050 590}%
\special{pa 4940 2020}%
\special{fp}%
%
\special{pn 8}%
\special{sh 1}%
\special{ar 4050 2010 10 10 0  6.28318530717959E+0000}%
\special{sh 1}%
\special{ar 4050 2010 10 10 0  6.28318530717959E+0000}%
%
\special{pn 8}%
\special{sh 1}%
\special{ar 4040 600 10 10 0  6.28318530717959E+0000}%
\special{sh 1}%
\special{ar 4040 600 10 10 0  6.28318530717959E+0000}%
%
\special{pn 8}%
\special{sh 1}%
\special{ar 3180 2010 10 10 0  6.28318530717959E+0000}%
\special{sh 1}%
\special{ar 3180 2010 10 10 0  6.28318530717959E+0000}%
%
\special{pn 8}%
\special{sh 1}%
\special{ar 4950 2010 10 10 0  6.28318530717959E+0000}%
\special{sh 1}%
\special{ar 4950 2010 10 10 0  6.28318530717959E+0000}%
%
\special{pn 8}%
\special{ar 1490 1520 1160 780  0.6823798 2.4489387}%
%
\special{pn 8}%
\special{sh 1}%
\special{ar 1500 2300 10 10 0  6.28318530717959E+0000}%
\special{sh 1}%
\special{ar 1500 2300 10 10 0  6.28318530717959E+0000}%
%
\special{pn 8}%
\special{sh 1}%
\special{ar 1950 2240 10 10 0  6.28318530717959E+0000}%
\special{sh 1}%
\special{ar 1950 2240 10 10 0  6.28318530717959E+0000}%
%
\special{pn 8}%
\special{sh 1}%
\special{ar 4520 2010 10 10 0  6.28318530717959E+0000}%
\special{sh 1}%
\special{ar 4520 2010 10 10 0  6.28318530717959E+0000}%
\put(15.0000,-4.2000){\makebox(0,0)[lb]{$a_{m+1}$}}%
\put(41.4000,-4.2000){\makebox(0,0)[lb]{$a_{m+1}'$}}%
\put(3.1000,-20.1000){\makebox(0,0)[lb]{$z_0$}}%
\put(25.1000,-20.1000){\makebox(0,0)[lb]{$z_1$}}%
\put(29.2000,-20.1000){\makebox(0,0)[lb]{$z_0'$}}%
\put(51.1000,-20.1000){\makebox(0,0)[lb]{$z_1'$}}%
\put(14.9000,-18.6000){\makebox(0,0)[lb]{$b$}}%
\put(14.2000,-25.1000){\makebox(0,0)[lb]{$b_m$}}%
\put(19.9000,-24.2000){\makebox(0,0)[lb]{$z$}}%
\put(39.6000,-22.3000){\makebox(0,0)[lb]{$b_m'$}}%
\put(44.4000,-22.5000){\makebox(0,0)[lb]{$z'$}}%
\put(8.3000,-23.7000){\makebox(0,0)[lb]{$c$}}%
\end{picture}%
\caption{$\triangle(z_0,z_1,a_{m+1})$ and
$\triangle(z_0',z_1',a_{m+1}')$}
\label{fig:smallerp}
  \end{figure}
  
  For any given two points $x_0',x_1' \in \triangle_1^{m+1}$, either
  segment $[a_{m+1}',x_i']$ extends to
  a segment $[a_{m+1}',\overline{x}_i']$ with $\overline{x}_i' \in \triangle_1^m$.
  Let $x_i := \varphi_{m+1}(x_i')$ and $\overline{x}_i :=
  \varphi_{m+1}(\overline{x}_i')$.  Since $\varphi_m$ is $1$-Lipschitz
  continuous, we have $d(\overline{x}_0,\overline{x}_1) \le
  d(\overline{x}_0',\overline{x}_1')$.  Comparing
  $\triangle(\overline{x}_0,\overline{x}_1,a_{m+1})$ and
  $\triangle(\overline{x}_0',\overline{x}_1',a_{m+1}')$, the
  $\CAT(1)$-inequality leads to $d(x_0,x_1) \le d(x_0',x_1')$.  Thus,
  $\varphi_{m+1}$ is $1$-Lipschitz continuous.  It remains to prove
  that
  \begin{equation}
    d(x_0,x_1) > d(x_0',x_1') - O(\epsilon^{1/8}).
    \label{eqn:smaller1}
  \end{equation}
  By Lemma \ref{lem:discurve}, for any point $w \in
  [\overline{x}_0,\overline{x}_1]$ there exists a point $z \in
  \varphi_m([\overline{x}_0',\overline{x}_1'])$ with $d(z,w) <
  2\epsilon$.  This and \eqref{eqn:smallera} imply
  \begin{equation}
    d(a_{m+1},w) > \pi/2 - O(\epsilon^{1/2}).
    \label{eqn:smallerb}
  \end{equation}
  Consider the small geodesic triangles
  $\triangle(\overline{x}_0,\overline{x}_1,a_{m+1})$ in $Y$ and
  $\triangle(\overline{x}_0',\overline{x}_1',a_{m+1}')$ in a unit
  $2$-sphere in $\triangle_1^{m+1}$.  By $\dis \varphi_m < \epsilon$,
  \eqref{eqn:smallerb}, and applying Lemma \ref{lem:ruled} to their
  triangles, we obtain \eqref{eqn:smaller1}.  This completes the proof
  of Lemma \ref{lem:induction}.
\end{proof}

To prove Theorem \ref{thm:smaller}, we set $n := \dim_C Y + 1 <
\infty$ and suppose that $\rad Y = \pi/2$.  Note that $\rho(y) =
\pi/2$ holds for any $y \in Y$.  Let $\epsilon$ be a positive number
with $\epsilon \ll \delta_n$.  Take a point $a_0 \in Y$.  There exists
a point $a_1 \in Y$ with $d(a_1,a_0) > \pi/2 - \epsilon$.  Let
$\varphi_1 : \triangle_1^1 \to [a_0,a_1]$ be the linear bijective map.
Since $\pi/2 - \epsilon < d(a_0,a_1) \le \pi/2$, this is a
$1$-Lipschitz continuous map $\triangle_1^1 \to Y$ with $\dis
\varphi_1 < \epsilon$.  By Lemma \ref{lem:induction}, we inductively
have $1$-Lipschitz continuous maps $\varphi_m \colon \triangle_1^m \to
Y$, $m = 1, 2, \dots, n$, such that $\dis \varphi_m <
O(\epsilon^{1/8^m})$.  Since $\dim \varphi_n(\triangle_1^n) \le \dim_C
Y = n-1$, Lemma \ref{lem:sperner} implies that $\varphi_n$ collapses
$\partial \triangle_1^n$.  Hence, there exist $n+1$ points $y_i' \in
F_i$, $i = 0, 1, \dots, n$, that are all mapped by $\varphi_n$ to a
common point of $Y$, where $F_i \subset \triangle_1^n$ is the opposite
face to $a_i'$.  We set
$$
\alpha_n := \inf \{ \, \max_{i, j} d(x_i',x_j') \mid \ x_i' \in F_i, \
x_i' \neq x_j' \, \text{for any} \, i \neq j \, \} > 0.
$$
Then for some $i_0 \neq j_0$ we have
$\alpha_n \le d(y_{i_0}',y_{j_0}') \le \dis \varphi_n < O(\epsilon^{1/8^n})$,
which is a contradiction if $\epsilon$ is sufficiently small.
Therefore we obtain $\rad Y < \pi/2$.

Since $Y$ is complete, $C(Y)$ consists of a single point
(cf.~Proposition 3.1 in \cite{langsch2}).
This completes the proof of Theorem \ref{thm:smaller}.
\qed

\medskip

As mentioned in Section 1,
we have proved that for $Y$ as in Theorem \ref{thm:smaller}
there exists a constant $\delta > 0$ depending only on $\dim_C Y < \infty$
such that $\rad Y \le \pi/2 - \delta$.

\begin{rem}\label{rem:slim}
  For $A \subset Y$, we denote by $C_Y(A)$ the set of all points where
  the function $Y \ni x \mapsto \sup_{y \in A} d(x,y) \in [0,\infty]$
  attains the infimum.  For an arbitrary subset $A$ of a
  $\CAT(1)$-space $Y$ with $\diam A \le \pi/2$, we have $\diam A =
  \diam B$ for the closure $B$ of the convex hull of $A$ (cf.~Lemma
  4.1 in \cite{langsch1}).  By applying Theorem \ref{thm:smaller} to
  $B$, we obtain the following generalization.  Let $Y$ be a complete
  $\CAT(1)$-space of $\dim_C Y < \infty$, and $A \subset Y$ a subset
  of $\diam A \le \pi/2$.  Then $\inf_{x \in Y} \sup_{y \in A} d(x,y)
  < \pi/2$, and $C_Y(A)$ consists of a single point.
\end{rem}

\subsection
{Proof of Theorem \ref{thm:dcenter}}

Let $f$ be a parabolic isometry of a proper $\CAT(0)$-space $X$ and
let $B := X_f(\infty)$.
It follows from 
Theorem \ref{thm:radfix} and Lemma \ref{lem:monoset}(3) that
$B$ is a closed, $\pi$-convex subset
of $(X(\infty),Td)$ with $\rad B \le \pi/2$.
Hence, $B$ itself is a complete $\CAT(1)$-space. 

First, we verify that $C(B)$ is non-empty.  Let $\rho : B \to
\mathbb{R}$ be the function defined by $\rho(x) := \sup_{y \in B}
d(x,y)$.  There exists a sequence $\{x_i\}$ in $B$ with $\rho(x_i) \to
\rad B$ as $i \to \infty$.  Since $X(\infty)$ is compact with respect
to the sphere topology, some subsequence of $\{x_i\}$ converges to a
point $x$.  We have $x \in B$ because $B$ is closed.  By the lower
semi-continuity of Tits distances, we have $\rho(x) \le \rad B$.
Thus, $\rho(x) = \rad B$ and $x \in C(B)$.

By the convexity of $B$, $C(B)$ is a closed, convex subset of
$(X(\infty),Td)$ with the property that $\diam C(B) \le \rad B \le
\pi/2$.  By setting $Y := C(B)$, it is a complete $\CAT(1)$-space of
$\diam Y \le \pi/2$.  By Proposition \ref{prop:cptdimension}, $\dim_C
(X(\infty),Td)$ is finite.  Therefore, by Theorem \ref{thm:smaller} we
have $\rad Y < \pi/2$, and $C(Y)$ consists of a single point.
Moreover, the second half follows from the uniqueness of the point and
its property.  This completes the proof of Theorem \ref{thm:dcenter}.
\qed

\medskip

Theorem \ref{thm:smaller} and the proof of Theorem \ref{thm:dcenter}
imply the following:

\begin{prop}\label{prop:ccdcenter}
  Let $Y$ be a compact $\CAT(1)$-space of $\dim Y < \infty$ and $\rad
  Y \le \pi/2$.  Then $C^2(Y)$ consists of a single point.
\end{prop}

\begin{rem}\label{rem:admissible}
  Let $X$ be a complete $\CAT(0)$-space and $G$ a subgroup of the
  isometry group of $X$.  Set $X_G(\infty) := \bigcap \{ X_g(\infty) |
  g \in G \}$.  We say that $G$ is \emph{admissible} if $X_G(\infty)
  \neq \emptyset$ and $\rad X_G(\infty) \le \pi/2$.  It follows from
  Theorem \ref{thm:dcenter} that if $G$ is an abelian group containing
  a parabolic element, then $G$ is admissible, provided $X$ is as in
  Theorem \ref{thm:dcenter}.  This is an extension of Proposition
  4.4.2 of \cite{eberlein}.  Similarly, we can obtain some extensions
  of Propositions 4.4.3, 4.4.4, and Corollary 4.4.5 of \cite{eberlein}
  for $\CAT(0)$-spaces.  Proposition 4.4.6 of \cite{eberlein} can be
  also extended by using the flat torus theorem for $\CAT(0)$-spaces
  (cf.~Theorem II.7.1 of \cite{bridhaef}).
\end{rem}

\section
{Example of a symmetric space}
 
In this section we discuss the symmetric space 
$SL(3,{\Bbb R})/SO(3,{\Bbb R})$ 
in detail as an example for Theorems \ref{thm:radfix} and \ref{thm:dcenter}.
A good reference for standard facts we use here is 
II.10 in \cite{bridhaef}. We would like to thank M.~Bestvina
for suggesting this example, and also K.~Wartman for 
useful discussions and informations.

\subsection{Manifolds $P(n,\mathbb{R})$ and $P(n,{\Bbb R})_1$}
Let $P(n,\mathbb{R})$ denote 
the space consisting of all positive definite, symmetric $(n \times n)$-matrices
with real coefficients.
Naturally, $P(n,\mathbb{R})$ is a differentiable manifold
of dimension $n(n+1)/2$.
The tangent space $T_pP(n,{\Bbb R})$
at a point $p$ is naturally isomorphic (via translation)
to the space of all symmetric $(n \times n)$-
matrices, $S(n,{\Bbb R})$.
The inner product $(u,v)_p = \tr (p^{-1}up^{-1}v)$ on 
$T_pP(n,{\Bbb R}) \simeq S(n,{\Bbb R})$
defines a Riemannian metric on $P(n,{\Bbb R})$, 
where $\tr u$ is the trace of a matrix $u$.
$P(n,\mathbb{R})$ is a simply connected, complete, non-positively 
curved Riemannian manifold, so that it is a proper $\CAT(0)$-space.

Let $P(n,\mathbb{R})_1 \subset P(n,\mathbb{R})$  be the subset 
of matrices with determinant 1.
$P(n,\mathbb{R})_1$ is a totally geodesic submanifold, 
whose tangent space at $p$ is the subspace
in $S(n,{\Bbb R})$ of matrices with trace 0.
$P(n,\mathbb{R})$ is a simply connected, complete, non-positively 
curved Riemannian manifold of dimension 
$n(n+1)/2-1$, so that it is a proper $\CAT(0)$-space as well.

$SL(n,\mathbb{R})$ acts on $P(n,\mathbb{R})$
by isometries according to the rule
$$f(p) := fp\,^t\!f, \,  (p \in P(n,\mathbb{R}), f \in SL(n,\mathbb{R})),$$
where $^t\!f$ is the transpose of $f$.
The right hand side of the definition is by the multiplication of matrices.
We may write $f \cdot p$ instead of $f(p)$.
$P(n,\mathbb{R})_1$ is invariant by this action, and 
the action is transitive on this submanifold.
Let $e$ be the identity matrix. 
The stabilizer of $e$ is $SO(n)$, so that 
$P(n,\mathbb{R})_1$ is identified as $SL(n,\mathbb{R})/SO(n)$.

\subsection
{Geometry of $P(3,\mathbb{R})_1$ and Tits boundary}
We collect some standard facts on $P(3,\mathbb{R})_1$
 from II.10 in \cite{bridhaef}. 
Most of them are true for all $P(n,\mathbb{R})_1, n \ge 3$
with appropriate change. 
Put $X: = P(3,\mathbb{R})_1$.
$X$ is a $5$-dimensional,
irreducible symmetric space of non-compact type of rank $2$,
which is a proper $\CAT(0)$-space.

Let us denote the Tits boundary $(X(\infty),Td)$ by
$X(\infty)$ for simplicity.
$X(\infty)$ is 
a ``thick spherical building" of dimension $1$ such that 
each apartment is isometric to $\mathbb{S}^1(1)$ and
each  Weyl chamber at infinity is an edge of length $\pi/3$.
$\diam X(\infty) = \pi$.
Since $X(\infty)$ is a spherical building,  
any two Weyl chambers at infinity are contained in at least one apartment.

The action of 
$SL(3,\mathbb{R})$ induced on $X(\infty)$ is 
by simplicial isometries.
It is transitive on
pairs $(A,c)$, where $A$ is an apartment, and $c \subset A$ is 
a Weyl chamber at infinity.
A Weyl chamber is a fundamental domain for the action.
(cf. II.10.71,75,76,77 in \cite{bridhaef}).
Therefore there are two orbits in the vertices of
$X(\infty)$ by the group action, so that 
$X(\infty)$ is a bi-partite
graph. It follows that any loop in $X(\infty)$ consists
of an even number of edges.

The isometry group of $X$, $I(X)$, has 
two connected components, and the one 
which contains the identity map, $I_0(X)$,  
 is $SL(3,\mathbb{R})$.
Let $\sigma$ be the involution of $X$ at $e$, which 
is an orientation reversing isometry.
It is given by $\sigma(f)=\, ^t\!f^{-1}$.
$I(X)=I_0(X)  \cup \sigma I_0(X)$.

Let $f$ be an isometry of $X$.
$\Min(f)$ denotes the set of all points in $X$ at which
the displacement function $d_f$ of $f$ attains its infimum
$|f| := \inf_{p \in X} d_f(p)$, 
which is the translation length.
If $f$ is elliptic, then $\Min(f)$ coincides with
the fixed point set $\Fix(f)$ of $f$ in $X$.
If $f$ is hyperbolic, then
the axes of $f$ are parallel to each other, and the union of their images
is $\Min(f)$.
If $f$ is parabolic, then $\Min(f) = \emptyset$. 
$f$ is said to be {\it semi-simple}
if $f$ is elliptic or hyperbolic.

In this section we calculate those geometric characters
of  $f \in SL(3,\mathbb{R})$.

\subsection{Real Jordan forms}\label{jordan}
It is known that $f \in SL(3,\mathbb{R})$ is semi-simple as an isometry of $X$
if and only if it is semi-simple as a matrix, i.e., diagonalizable
in $GL(3,{\Bbb C})$.
(cf. II.10.61 in \cite{bridhaef}).

Calculation of  $X_f(\infty)$ and $\Min(f)$ of $f \in SL(3,\mathbb{R})$
is mostly by linear algebra.
Each $f \in SL(3,\mathbb{R})$ 
is conjugate to $g$ 
in $SL(3,\mathbb{R})$ such that $g$ is 
one (and only one) of the following list.
$g$  is 
a real Jordan form of $f$.
The symbol $\diag(a,b,c)$ is for the $(3 \times 3)$-diagonal 
matrix with entries $a,b,c$.


Since $f$ and $g$ are conjugate in $I(X)$, 
$f$ is elliptic, hyperbolic, or parabolic
if and only if so is $g$, respectively.
If $h \in I(X)$ is a conjugating element, i.e., 
$hfh^{-1}=g$, then 
$
X_f(\infty) = h \cdot X_g(\infty)$, 
$\Min(f) = h \cdot \Min(g),
$
and $|f| = |g|$.
We discuss $g$ instead of $f$.


\medskip
{\bf List of real Jordan forms in $SL(3,{\Bbb R})$}.
\begin{enumerate}
\item
$
\begin{pmatrix}
1 & 0 & 0 \\
0 & 1 & 1 \\
0 & 0 & 1
\end{pmatrix}.
$
\item
$
\begin{pmatrix}
1/a^2 & 0 & 0 \\
0 & a & 1 \\
0 & 0 & a
\end{pmatrix}, 
$
where  
$0, 1 \not = a \in {\Bbb R}$.
\item
$
\begin{pmatrix}
1 & 1 & 0 \\
0 & 1 & 1 \\
0 & 0 & 1
\end{pmatrix}.
$
\item
$
\begin{pmatrix}
a & b & 0 \\
-b & a & 0 \\
0 & 0 & 1/(a^2+b^2)
\end{pmatrix}
\label{eqn: cpx}
$
where $a, b \in \mathbb{R}$ with $a^2+b^2 \not =0$ and $b \neq 0$.

This one  is conjugate to $\diag(a+ib,a-ib,1/(a^2+b^2))$
by an element in $GL(3,{\Bbb C})$.

\item
$\diag(a,b,c)$
such that $a,b,c \in {\Bbb R}, abc=1, a\not=b \not = c  \not =a$.
\item
$\diag(a,a,1/a^2), a \in {\Bbb R}, a \not = 0,1$.
\item
$e.$
\end{enumerate}

\subsection{Flat and Weyl chambers}\label{ss.flat}

Consider the following linear subspace 
in $T_eX$.
$$\mathfrak{a}_0 := \{ u \mid u = \diag (u_1,u_2,u_3), \, \tr u = 0 \} \subset T_eX.$$
Let 
$$
F_0 := \{ \ \exp (u) \mid  u 
\in \mathfrak{a}_0 \, \} \subset X.
$$
$F_0$ is a flat plane in $X$
and
$A_0 := F_0(\infty)$ is an apartment in $X(\infty)$.


For $x \in A_0$,
let $\gamma_{ex}$ be the 
geodesic in $F_0$ from $e$ to $x$.
$\gamma_{ex}$
is $\exp(t u(x)), t \ge 0$
for some $u(x) \in \mathfrak{a}_0$.
The tangent vector at $e$, $u(x)$, is uniquely determined by $x$ upto
scaling by a positive number, so that 
let us denote the one of unit length by 
$u(x)$, 
which we write as
$\diag(u_1(x),u_2(x),u_3(x))$.

\begin{figure}[htbp]
  \centering
\unitlength 0.1in
\begin{picture}(19.10,21.50)(11.30,-24.00)
%
\special{pn 13}%
\special{pa 2110 470}%
\special{pa 2940 950}%
\special{pa 2939 1909}%
\special{pa 2109 2388}%
\special{pa 1279 1908}%
\special{pa 1279 949}%
\special{pa 2110 470}%
\special{fp}%
%
\special{pn 8}%
\special{pa 1270 1890}%
\special{pa 2940 950}%
\special{fp}%
\put(19.9000,-14.0000){\makebox(0,0)[lb]{$\gamma_0$}}%
\put(20.5000,-4.2000){\makebox(0,0)[lb]{$v_1$}}%
\put(30.4000,-9.0000){\makebox(0,0)[lb]{$v_2$}}%
\put(30.3000,-19.9000){\makebox(0,0)[lb]{$v_3$}}%
\put(20.7000,-25.7000){\makebox(0,0)[lb]{$v_4$}}%
\put(11.6000,-20.6000){\makebox(0,0)[lb]{$v_5$}}%
\put(11.3000,-9.2000){\makebox(0,0)[lb]{$v_6$}}%
\put(24.4000,-21.2000){\makebox(0,0)[lb]{$c_3$}}%
\put(16.8000,-20.9000){\makebox(0,0)[lb]{$c_4$}}%
\put(13.4000,-14.7000){\makebox(0,0)[lb]{$c_5$}}%
\put(16.6000,-9.1000){\makebox(0,0)[lb]{$c_6$}}%
\put(23.3000,-8.7000){\makebox(0,0)[lb]{$c_1$}}%
\put(27.9000,-15.5000){\makebox(0,0)[lb]{$c_2$}}%
\end{picture}%
  \caption{6-gon}
  \label{fig:6gon}
\end{figure}

$A_0$ is a $6$-gon as a building
with the following Weyl chambers (see Figure \ref{fig:6gon}).
$\{ u_i(x) \ge u_j(x) \ge u_k(x) \}$
means the set 
 $\{ x \in A_0 \mid u_i(x) \ge u_j(x) \ge u_k(x) \}$.
\begin{align*}
c_1 &:= \{ \, u_1(x) \ge u_2(x) \ge u_3(x) \, \}, \quad
c_2 := \{ \, u_2(x) \ge u_1(x) \ge u_3(x) \, \}, \\
c_3 &:= \{ \, u_2(x) \ge u_3(x) \ge u_1(x) \, \}, \quad
c_4 := \{ \, u_3(x) \ge u_2(x) \ge u_1(x) \, \}, \\
c_5 &:= \{ \, u_3(x) \ge u_1(x) \ge u_2(x) \, \}, \quad
c_6 := \{ \, u_1(x) \ge u_3(x) \ge u_2(x) \, \}.
\end{align*}
$A_0 = c_1 \cup \cdots \cup c_6$.
Define $v_i := c_{i-1} \cap c_i$, 
$i =  1, 2, \dots, 6$, where $c_0=c_6$.
They are the vertices of $X(\infty)$ in $A_0$.
We may write $c_i=[v_i,v_{i+1}], 1 \le i \le 6$, where $v_7=v_1$.

A bi-infinite geodesic, or simply line, in $X$ is always contained in some 
flat plane because $X$ is a symmetric space.
If  a line  is contained in a unique flat, then it is called 
regular (cf. 10.46 \cite{bridhaef}), and otherwise it is called singular.
There are three singular bi-infinite geodesics (without orientation)
on $F_0$, which are $\gamma_{ev_1} \cup \gamma_{ev_4},
\gamma_{ev_2} \cup \gamma_{ev_5}, \gamma_{ev_3} \cup \gamma_{ev_6}$.

Set $w_2 :=u(v_2)= (1/\sqrt{6}) \diag(1,1,-2)$.
This is a unit vector at $e$ tangent to $F_0$, 
pointing the vertex $v_2$ at infinity.
Define a line in $F_0$ by 
$$\gamma_0(t) := \exp (tw_2), t \in {\Bbb R}.$$
$\gamma_0$ is a line 
through $e$ such that $\gamma_0(\infty)=v_2, 
\gamma_0(-\infty)=v_5$.
As a set, $\gamma_0=\{ \diag(s,s,1/s^2)|0 < s \in {\Bbb R} \}
= \gamma_{ev_2} \cup \gamma_{ev_5}$.
$\gamma_0$ is a singular geodesic.

For a line $\gamma$ in $X$, let 
$P(\gamma)$ denote the union of all lines in $X$ 
parallel to $\gamma$. This is a convex subset in $X$, 
so that let $P(\gamma)(\infty) \subset X(\infty)$ denote the set 
of points at infinity of $P(\gamma)$.

Denote by $\mathcal{F}_0$ the set of all flat planes in $X$
containing $\gamma_0$.
Then, 
$P(\gamma_0) =  \bigcup \{  F \mid F \in \mathcal{F}_0  \}$.
$P(\gamma_0)$
is a totally geodesic, $3$-dimensional submanifold, 
which is naturally 
isometric to $P(2,\mathbb{R})_1 \times \mathbb{R}$
(cf.~Proposition II.10.67 in \cite{bridhaef}).
$P(2,\mathbb{R})_1$ is isometric to
$\mathbb{H}^2$ up to a scaling factor.
We note that 
$P(\gamma_0)(\infty) = \bigcup \{  F(\infty) \mid 
F \in \mathcal{F}_0  \}$.

\subsection{Theorem}

\begin{thm}\label{sl3}
Suppose $g \in SL(3,{\Bbb R})$ is one 
on the list in the subsection \ref{jordan}.
Then we have the following in the order on the list:
\begin{enumerate}
\item
$g$ is parabolic and $X_g(\infty)$ is the union of all edges incident to $c_2$.
$X_g(\infty)$ is not compact in $(X(\infty),Td)$, with 
uncountably many edges.
$|g|=0$.
\item
$g$ is parabolic and $X_g(\infty)=c_1\cup c_2\cup c_3$.
$|g|=2\sqrt{6}\log|a|$.

\item
$g$ is parabolic and $X_g(\infty)=c_1$.
$|g|=0$.

\item
$g$ is semi-simple, and 
\begin{enumerate}
\item
$X_g(\infty)=\{ v_2, v_5 \}$.
\item
If $a^2+b^2=1$, then 
$g$ is elliptic, 
and $\Fix(g) = \gamma_0$.
\item
If $a^2+b^2 \not =1$, then 
$g$ is hyperbolic and 
$|g|=\sqrt{6} \log(a^2+b^2).$
$\Min(g)=\gamma_0$.
\end{enumerate}

\item
$g$ is hyperbolic, and 
$|g|=2\sqrt{(\log|a|)^2+(\log|b|)^2+(\log|c|)^2}$.
\begin{enumerate}
\item
$X_g(\infty)=A_0$.
\item
$\Min(g)=F_0$.
\end{enumerate}

\item
$g$ is hyperbolic, and 
$|g|=2\sqrt{6}\log|a|$.
\begin{enumerate}
\item
$X_g(\infty)=P(\gamma_0)(\infty)$.
\item
$\Min(g)=P(\gamma_0)$.
\end{enumerate}

\item
$g$ is the identity map, so that elliptic, 
with $\Fix(g)=X$ and $X_g(\infty)=X(\infty)$.

\end{enumerate}

\end{thm}

\subsection{Stabilizers}
The analysis of the 
stabilizing subgroup in $SL(3,{\Bbb R})$
of a point $v \in X(\infty)$ is important for the proof of the 
theorem.
We quote  Proposition II.10.64 in \cite{bridhaef}
in the following form.

\begin{lem}\label{fix}
Let $g = (g_{ij}) \in SL(3,\mathbb{R})$, and
$x \in A_0$.
Then $g(x) = x$ if and only if
$g_{ij}e^{-t(u_i(x)-u_j(x))}$ 
converges as $t \to \infty$ 
 for all $i, j$.
\end{lem}

This implies the following.
\begin{prop}\label{stab}
Let $G_i$ be the 
 subgroups of $SL(3,{\Bbb R})$ stabilizing $v_i$.
Then, 
$$ 
G_1=
\left\{
\begin{pmatrix}
* & * & * \\
0 & * & * \\
0 & * & *
\end{pmatrix}
 \right\}
,
G_2=
\left\{
\begin{pmatrix}
* & * & * \\
* & * & * \\
0 & 0 & *
\end{pmatrix}
\right\}
,
G_3=
\left\{
\begin{pmatrix}
* & 0 & * \\
* & * & * \\
* & 0 & *
\end{pmatrix}
\right\}
,
$$
$$
G_4=
\left\{
\begin{pmatrix}
* & 0 & 0 \\
* & * & * \\
* & * & *
\end{pmatrix}
\right\}
,
G_5=
\left\{
\begin{pmatrix}
* & * & 0 \\
* & * & 0 \\
* & * & *
\end{pmatrix}
\right\}
,
G_6=
\left\{
\begin{pmatrix}
* & * & * \\
0 & * & 0 \\
* & * & *
\end{pmatrix}
\right\}
,
$$
where $* \in {\Bbb R}$.
\end{prop}

Let $H_1$ be the following  subgroup, 
parameterized by $t \in {\Bbb R}$, 
which fixes edges $c_1,c_2,c_3$, pointwise.

$$
H_1=
\left\{
\begin{pmatrix}
1 & 0 & 0 \\
0 & 1 & t \\
0 & 0 & 1
\end{pmatrix} \Biggl|\Biggr.\ t \in {\Bbb R}
\right\}.
$$

$H_1$ fixes $v_1 \in X(\infty)$. 
$H_1$ acts transitively 
on the set of all edges 
incident to $v_1$ other than $c_1$.
To see it, 
consider the following subgroup in $SL(3,{\Bbb R})$ 
containing $H_1$.
$$J_1=
\left\{
\begin{pmatrix}
* & 0 & 0 \\
0 & * & * \\
0 & 0 & *
\end{pmatrix} \Biggl|\Biggr.\ * \in {\Bbb R}
\right\}
.
$$
For a given edge $c \not = c_1$, incident to $v_1$, 
we will find $h \in H_1$ with $h(c)=c_6$.
Take an apartment, $A$, containing
$c$ and $c_3$. Then it automatically contains
$c_1,c_2$ as well.
Recall that 
$SL(3,{\Bbb R})$ acts transitively 
on the set of pairs of an apartment, $A'$, in $X(\infty)$
and a Weyl chamber, $c'$, in $A'$, $(A',c')$.
Take $j \in SL(3, {\Bbb R})$ which 
maps $(A,c)$ to $(A_0,c_6)$. Clear that $j \in J_1$
since it fixes $c_1,c_2,c_3$.
Let
$j=
\begin{pmatrix}
p & 0 & 0 \\
0 & q & s \\
0 & 0 & r
\end{pmatrix}
, pqr=1$.
Take 
$k=
\begin{pmatrix}
1/p & 0 & 0 \\
0 & 1/q & 0 \\
0 & 0 & 1/r
\end{pmatrix}
\in SL(3,{\Bbb R})$.
Then 
$kj=
\begin{pmatrix}
1 & 0 & 0 \\
0 & 1 & s/q \\
0 & 0 & 1
\end{pmatrix}
=h 
\in H_1$.
We have $k(c_i)=c_i$ for all $i$, 
so it follows from $j(c)=c_6$
that $h(c)=kj(c)=k(c_6)=c_6$.

Also, $H_1$ acts transitively on 
 the set of all edges incident
at $v_4$ other than $c_3$. Each of the two sets 
is parameterized by $t$ by the action.

Off-diagonal entries except for the $(2,3)$-entries
of matrices in $H_1$ are $0$.
Since there are 6 off-diagonal entries in $(3 \times 3)$-matrices, 
we consider 5 other similar subgroups, $H_2,H_3,H_4,H_5,H_6$, 
which we define later.

\subsection{Proof}
We discuss each case in the order and prove Theorem \ref{sl3}.
Case 7 is trivial.

\medskip
{\bf Case 1}.
%
$g$ is parabolic since it is not diagonalizable
as a matrix.
By Proposition \ref{stab}, 
$X_g(\infty) \cap A_0 =c_1\cup c_2\cup c_3$.
To see that any edge, $c \not =c_2, c_3$, incident to $v_3$ is fixed by $g$, take 
a (unique) element $h \in H_6$ such that $h(c)=c_3$, where
$$
H_6=
\left\{
\begin{pmatrix}
1 & 0 & t \\
0 & 1 & 0 \\
0 & 0 & 1
\end{pmatrix} \Biggl|\Biggr.\ t \in {\Bbb R}
\right\}
.
$$
Since $g(c_3)=c_3$, $h^{-1}gh(c)=c$. Then
it follows from $hg=gh$ that $g(c)=c$, pointwise.

To see any edge $c \not =c_2$, incident to $v_2$ is fixed by $g$, take 
a (unique) element $h \in H_2$ such that $h(c)=c_1$, where
$$
H_2=
\left\{
\begin{pmatrix}
1 & 0 & 0 \\
t & 1 & 0 \\
0 & 0 & 1
\end{pmatrix} \Biggl|\Biggr.\ t \in {\Bbb R}
\right\}.
$$
As before, we have $hg=gh$, therefore $g(c)=c$, pointwise.
We know $X_g(\infty)$ has no more edges because
it is connected and its diameter is at most $\pi$.

The edges in $X_g(\infty)$ other than 
$c_1,c_2,c_3$ are those which are parametrized by 
$H_2$ and the others which are parametrized by $H_6$, so that uncountable.
It is not compact because the mid points of the edges
are at least $\pi/3$ apart from each other.

Let us prove $|g|=0$ by computation.
To deal with Case 3 at one time, suppose
$g=\begin{pmatrix}
1 & k & 0 \\
0 & 1 & 1 \\
0 & 0 & 1
\end{pmatrix}.$
It suffices for us to show that there exists a geodesic, $\gamma(t)$, 
such that $\lim_{t \to \infty} d(g(\gamma(t)),\gamma(t))=0$.
We use the notations from the subsection \ref{ss.flat}.
For $x \in A_0$, set $\gamma(t)=\gamma_{ex}(t)=
\exp(tu)$, where $u$ is the diagonal matrix $u(x)=\diag(u_1,u_2,u_3)$.
For simplicity, we write the result of the action 
by a group element $g$ on a point $p$ as $g \cdot p$, instead of $g(p)$, 
in this discussion.
Then,
$$d(g \cdot \exp(tu),\exp(tu))=
d\left(\exp\left(-\frac{t}{2}u\right) \cdot g \cdot \exp(tu),
\exp\left(-\frac{t}{2}u\right) \cdot \exp(tu)\right)$$
$$=d\left(\exp\left(-\frac{t}{2}u\right)g\exp\left(\frac{t}{2}u\right) \cdot e,
e\right),$$
because $\exp\left(\frac{t}{2}u\right) \cdot e=\exp(tu)$.
By the computation of matrix multiplications,
$\exp\left(-\frac{t}{2}u\right)g\exp\left(\frac{t}{2}u\right)$
is
$$\begin{pmatrix}
1  & k\exp(t(u_2-u_1)/2) & 0 \\
0  & 1  & \exp(t(u_3-u_2)/2) \\
0 & 0  & 1
\end{pmatrix}.$$

For $g$ in Case 1, we have $k=0$, so that 
if $u_2 > u_3$, then as $t \to \infty$, this matrix tends to $e$, 
which means that $d(g \cdot \gamma(t), \gamma(t)) \to 0$.
We got $|g|=0$. We remark that $u_2 > u_3$
is satisfied for $x \in A_0$ 
if and only if $x \in (c_1 \cup c_2 \cup c_3) \backslash (v_1 \cup v_4)$.

For $g$ in Case 3, we have $k=1$.
So, if $u_1 > u_2 > u_3$, then $d(g \cdot \gamma(t), \gamma(t)) \to 0$, 
which shows that $|g|=0$.
The condition $u_1 > u_2 > u_3$ holds for $x \in A_0$ if and only if
$x \in c_1 \backslash (v_1 \cup v_2)$.


\medskip
{\bf Case 2}.
%
%
As in Case 1, $g$ is parabolic and
 $X_g(\infty) \cap A_0 =c_1 \cup c_2 \cup c_3$.
To see there are not more edges than those in $X_g(\infty)$, suppose 
there was an edge, $c \not = c_2, c_3$, incident to $v_3$
with $g(c)=c$. Take, as before, $h \in H_6$
such that $h(c)=c_3$.
$$
h=
\begin{pmatrix}
1 & 0 & t \\
0 & 1 & 0 \\
0 & 0 & 1
\end{pmatrix}.
$$
Then, $ hgh^{-1}(c_3)=c_3$.
Since $c \not = c_3$, we have $t \not =0$, which 
is important to get a contradiction in this case.
By computation 
$$
hgh^{-1}=
\begin{pmatrix}
1/a^2 & 0 & t(a-1/a^2) \\
0 & a & 1 \\
0 & 0 & a
\end{pmatrix}.
$$
Since $t(a-1/a^2) \not =0$, 
$hgh^{-1}$ is not in $G_4$, so that 
does not fix $v_4 \in c_3$, a contradiction.

To see there is no edge $c \not =c_1, c_2$  at $v_2$
with $g(c)=c$, use $H_2$, as before.
If there was, take $h \in H_2$ with $h(c)=c_1$ such that 
$$
h=
\begin{pmatrix}
1 & 0 & 0 \\
t & 1 & 0 \\
0 & 0 & 1
\end{pmatrix}, t \not = 0.
$$
Then $hgh^{-1}(c_1)=c_1$, pointwise.
By computation, 
$$
hgh^{-1}=
\begin{pmatrix}
1/a^2 & 0 & 0\\
t(1/a^2-a)  & a & 1 \\
0 & 0 & a
\end{pmatrix}
$$
such that $t(1/a^2-a)  \not =0$, 
therefore $hgh^{-1} \not \in G_1$ does not fix $v_1$, which 
gives a contradiction since it is supposed to 
fix $c_1=[v_1,v_2]$.

We get the claim  because $X_g(\infty)$ is connected
and with diameter at most $\pi$.

We postpone the computation of $|g|$ until
the proof of Case 6. It is not a coincidence
that $|g|$ is the same number in Cases 2 and 6.

\medskip
{\bf Case 3}.
As in Case 1, $g$ is parabolic and 
 $X_g(\infty) \cap A_0 =c_1$.
To see this is all, suppose 
there was an edge, $c \not =c_1$, incident to $v_1$
with $g(c)=c$. Take $h \in H_1$
such that $h(c)=c_6$.
It follows that $ hgh^{-1}(c_6)=c_6$.
Let 
$$
h=
\begin{pmatrix}
1 & 0 & 0 \\
0 & 1 & t \\
0 & 0 & 1
\end{pmatrix}.
$$
By computation 
$$
hgh^{-1}=
\begin{pmatrix}
1 & 1 & -t \\
0 & 1 & 1 \\
0 & 0 & 1
\end{pmatrix},
$$
therefore 
$hgh^{-1} \not \in G_6$ does not fix $v_6 \in c_6$, a contradiction.

To see $g$ does not fix any edge incident to $v_2$ other than 
$c_1$, use
$H_5$ and do the same argument, where 
$$
H_5=
\left\{
\begin{pmatrix}
1 & t & 0 \\
0 & 1 & 0 \\
0 & 0 & 1
\end{pmatrix} \Biggl|\Biggr.\ t \in {\Bbb R}
\right\}
.
$$
We get the claim  since $X_g(\infty)$ is connected.

Since we already showed that $|g|=0$ in the proof of Case 1, 
we finish Case 3.

\medskip
{\bf Case 4}.
%
%
$g$ is semi-simple because it is diagonalizable
in $GL(3,{\Bbb C})$.

\noindent
(a). By Proposition \ref{stab}, $X_g(\infty) \cap A_0 =\{ v_2, v_5 \}$.
To see this is all, we first show 
that there is no edge
incident to $v_2$ in $X_g(\infty)$.
Suppose there was one, $c$.
We know that $c \not = c_1, c_2$. Take $h \in H_2$
such that $h(c)=c_1$. 
If $
h=
\begin{pmatrix}
1 & 0 & 0 \\
t & 1 & 0 \\
0 & 0 & 1
\end{pmatrix} 
$, then 
$$
hgh^{-1}=
\begin{pmatrix}
a-tb & b & 0 \\
-b(1+t^2) & a+tb & 0 \\
0 & 0 & 1/(a^2+b^2)
\end{pmatrix},
$$
which does not fix $v_1$ because
$-b(1+t^2) \not =0$.
But $hgh^{-1}(c_1)=hg(c)=h(c)=c_1$, 
so that it fixes $v_1 \in c_1$, a contradiction.
Similarly there is no edge in $X_g(\infty)$ incident to $v_5$.


To finish, suppose there was a vertex, $v$, in $X_g(\infty) \backslash A_0$.
Then $Td(v,v_2)=\pi$, because
if it was less than $\pi$, then the unique geodesic from $v$ to $v_2$
would have to be in $X_g(\infty)$, which is impossible since 
there is no edge incident to $v_2$ fixed by $g$.
By the same reason, $Td(v,v_5)=\pi$.
Consider a loop made of three geodesics: one from $v$ to $v_2$, 
one from $v_2$ to $v_5$ and one from $v_5$ to $v$. This loop
consists of $9$ edges, which is impossible because $X(\infty)$ is a 
bi-partite graph.



(b).
Recall that 
$\gamma_0 = \{ \diag(s,s,1/s^2)|0 <  s \in {\Bbb R} \}$.
Under the condition 
$a^2+b^2=1$, by computation, $g p\, ^t\!g=p, p \in P(3,{\Bbb R})_1$
if and only if $p=\diag(s,s,1/s^2), 0 < s \in {\Bbb R}$.
We get the claim.

(c).
By computation, $g(\gamma_0)=\gamma_0$. 
Since $e \in \gamma_0$ and $e \not =g(e)$, 
$g$ is hyperbolic, and $\gamma_0$ is an axis.
The translation length, $|g|$, is 
$d(e, g(e))=d(e,ge^t\!g)=d(e,g^t\!g)$, where
$$g^t\!g=
\diag(a^2+b^2, a^2+b^2, 1/(a^2+b^2)^2)
=
\exp[ \log(a^2+b^2) \diag(1,1,-2)].
$$
Since the norm $|\diag(1,1,-2)|_e=\sqrt{6}$, 
$|g|=\sqrt{6} \log(a^2+b^2)$.

We know that $\gamma_0 \subset \Min(g) \subset P(\gamma_0)$ because
an axis of $g$ is parallel to $\gamma_0$.
Since $g$ leaves $\gamma_0$ invariant, 
$g$ leaves $P(\gamma_0)$ invariant as well.
We remark that the action of $g$ is 
by a shift and a rotation about $\gamma_0$.
We define the following subgroup in $SL(3,{\Bbb R})$, 
which is in fact in $SO(3)$.
\begin{multline*}
R = \{ \ h = (h_{ij}) \in SO(3) \mid  
(h_{ij})_{1 \le i,j \le 2} \in SO(2), \ h_{33} = 1, \\
h_{13}=h_{31}=h_{23}=h_{32}=0 
\ \}.
\label{eqn: stabilizer}
\end{multline*}
If $h \in R$, $h$ fixes
$e, v_2,v_5$, so that $h$ fixes all points
on $\gamma_0$.
Therefore $h$ leaves $P(\gamma_0)$ invariant, 
and acts on ${\mathcal F}_0$.

\smallskip
\noindent
{\it Claim}.
The action of $R$ on ${\mathcal F}_0$ 
is transitive.
\smallskip

To see it, let $F \in {\mathcal F}_0$ be
a flat. Then there is an element $w \in T_eX$
such that $w$ and $w_2$ commute 
as matrices and the image by $\exp$ of 
the subspace spanned by $w,w_2$ in $T_eX$ is 
$F$.
The two commuting symmetric matrices $w,w_2$
are simultaneously diagonalizable by an element, $h$, 
in $SO(3)$. Moreover since $w_2$ is 
diagonal, one may assume that $h$ commutes with $w_2$.
By computation, this implies that $h$ is in $R$.
We claim that $h$ maps $F$ to $F_0$. Indeed, 
let $\gamma$ be the geodesic through $e$ defined
by $\gamma=\exp(s w), s \in {\Bbb R}$. It is 
in $F$.
Since 
$h(\gamma_0)=\gamma_0$, it suffices to show
$h(\gamma) \subset F_0$. Since $h \in SO(3)$, 
$$h(\gamma)=h \exp(sw) ^t\!h=
h \exp(sw) h^{-1}=\exp (s hwh^{-1}), $$
which is in $F_0$
because $hwh^{-1}$ is diagonal.
We got the claim.

\medskip
Suppose there was an axis of $g$, $\alpha$, which 
is not $\gamma_0$. Take the plane $F \in {\mathcal F}_0$
which contains $\alpha$. Such $F$ exists since $\alpha$ is
parallel to $\gamma_0$.
Take $h \in R$ with $h(F)=F_0$.
Since $h$ commutes with $g$, $h(\alpha)$ is an axis of $g$ 
as well. It implies that $F_0$ is invariant by $g$, so that 
$F_0(\infty) \subset X_g(\infty)$, which is impossible.
We got $\Min(g)=\gamma_0$.
We finished Case 4.

\medskip
We are left with $g$ which are diagonal. 
Suppose 
$e \not =g = \diag (a,b,c) \in SL(3,\mathbb{R})$.
$g$ is a semi-simple isometry and 
the flat $F_0$ is $g$-invariant, so that 
$A_0 \subset X_g(\infty)$.

Set
$$N=\sqrt{(\log|a|)^2+(\log|b|)^2+(\log|c|)^2}$$
and define
a unit length element $u_g \in T_eF_0$ by 
$$
u_g := \frac{1}{N}
\diag ( \log |a|, \log |b|, 
\log |c|).
$$
Let $\gamma_g$ be the bi-infinite geodesic  in $F_0$ through $e$ defined by
$$\gamma_g(t) := \exp (tu_g), t \in {\Bbb R}.$$
This is of unit speed.
Computation shows that $\gamma_g$ is $g$-invariant, therefore
it is an axis. 
$|g|=d(e,g(e))=d(e,g^t\!g)$.
$g^t\!g=\diag(a^2,b^2,c^2)=\exp (2Nu_g)=\gamma_g(2N)$. 
Since $\gamma_g$ has unit speed, 
$|g|=2N$.

There are two cases: $\gamma_g$ is regular (Case 5) 
or singular (Case 6).
We already know that $g$ is hyperbolic and calculated
$|g|$.

\medskip
{\bf Case 5}.
%
%
%
(a).
Since $F_0$ is invariant by $g$, $A_0=F_0(\infty) \subset X_g(\infty)$.
Let $v \in X_g(\infty)$.
Then there is a flat $F$ with $\gamma_g \subset F$
and $v \in F(\infty)$.
Indeed, if $\gamma$ is a bi-infinite geodesic through
$e$ with $\gamma(\infty)=v$, then since $g(v)=v$, 
$\gamma$ and $\gamma_g$ is on some flat.

Since $\gamma_g$ is a regular geodesic, 
it is contained in only one flat, so that 
$F=F_0$. 
We get $v \in F_0(\infty)=A_0$.

(b). Since $F_0$ is $g$-invariant, $F_0 \subset \Min(g)$.
$\Min(g)$ consists of axes of $g$.
Let $\gamma$ be an axis different from $\gamma_g$.
Then there is a flat strip between them, so that 
there is indeed a flat, $F$, containing both of them
because it is in a symmetric space.
Since $\gamma_g$ is regular, we have $F=F_0$, so that 
$\gamma$ is in $F_0$.

\medskip
{\bf Case 6}.
%
%
%
(a). As in Case 5, $A_0 \subset X_g(\infty)$.
Since  $g$ commutes with any element in $R$,  
$X_g(\infty)$ is $R$-invariant, so that 
$R \cdot A_0 \subset X_g(\infty)$.
$R \cdot A_0=P(\gamma_0)(\infty)$ implies that 
$P(\gamma_0)(\infty) \subset X_g(\infty)$.
To see the other inclusion, let $v \in  X_g(\infty)$.
Then there is a flat, $F$, such that $\gamma_0 \subset F$
and $v \in F(\infty)$ (cf. (a) in Case 5).
By definition, $F \in {\mathcal F}_0$, so that 
$F \subset P(\gamma_0)$. We get $v \in P(\gamma_0)(\infty)$.

(b). In this case, $\gamma_g=\gamma_0$.
$F_0$ is $g$-invariant, so that $F_0 \subset \Min(g)$.
Since $g$ commutes with any element in $R$, 
$\Min(g)$ is $R$-invariant, so that 
$R \cdot F_0 \subset \Min(g)$.
Because  $R \cdot F_0=P(\gamma_0)$, 
$P(\gamma_0) \subset \Min(g)$.
On the other hand, since $P(\gamma_0)$ is 
the union of all geodesics parallel to $\gamma_0$, 
$\Min(g) \subset P(\gamma_0)$, therefore 
$\Min(g) = P(\gamma_0)$.
Case 6 is done. 

To finish the proof, we show  $|g|=2\sqrt{6}\log|a|$
for $g$ in Case 2.
It is easy to see that $g$ is conjugate in $SL(3,\Bbb{R})$
to the matrix, 
$h= \begin{pmatrix}
a & 1 & 0 \\
0 & a & 0 \\
0 & 0 & 1/a^2 
\end{pmatrix},
$
so that it suffices to show $|h|=2\sqrt{6}\log|a|$.

Although $\gamma_0$ is not $h$-invariant, 
$h$ fixes $\gamma_0(\infty)=v_2, \gamma_0(-\infty)=v_5$
because $h \in G_2 \cap G_5$, so that 
$h$ leaves not only the subset $P(\gamma_0)$, 
but also its product structure $P(\gamma_0)={\Bbb H}^2 \times {\Bbb R}$ invariant.

The restriction of $h$ to $P=P(\gamma_0)$, $h|_P$, is also parabolic.
Since $P$ is convex in $X$ and $h$-invariant, $|h|_P| = |h|$,
so that we compute $|h|_P|$.
$h|_P$ acts on $P={\Bbb H}^2 \times {\Bbb R}$ by 
a product of isometries:
a parabolic isometry on ${\Bbb H}^2$, denoted by $h|_{{\Bbb H}^2}$, and 
a translation on ${\Bbb R}$, denoted by $h|_{{\Bbb R}}$.
Since $|h|_{{\Bbb H}^2}|=0$, we have $|h|_P|=|h|_{{\Bbb R}}|$.

Consider the following matrix in $SL(3,{\Bbb R})$,
$k=
\begin{pmatrix}
1 & -1/a & 0 \\
0 & 1 & 0 \\
0 & 0 & 1 
\end{pmatrix}.
$
This is also a parabolic isometry, which 
leaves $P(\gamma_0)$ invariant such that it acts
on it as a product of isometries on ${\Bbb H}^2$ and 
${\Bbb R}$. The action of $k$ on ${\Bbb R}$
is trivial since the $(3,3)$-entry of $k$ is $1$.
This is because one can show by computation 
that the geodesic from $e \in X$ to $k(e)$ is 
perpendicular to $\gamma_0$  at $e$, 
so that $k(e) \in {\Bbb H}^2$ in ${\Bbb H}^2 \times {\Bbb R}$.
Or, one may use the fact that  $P(\gamma_0)$  is the union of 
matrices of the following form in $P(3,{\Bbb R})_1$;
$\begin{pmatrix}
* & * & 0 \\
* & * & 0 \\
0 & 0 & * 
\end{pmatrix},$
where the set of top-left $(2 \times 2)$-matrices
corresponds to ${\Bbb H}^2$ and the $(3,3)$-entries, 
which are  positive 
numbers, (by taking $\log$) correspond  to  ${\Bbb R}$ in the product 
decomposition
$P(\gamma_0) = {\Bbb H}^2 \times {\Bbb R}$.
By the definition of the action, 
$k$ acts trivially on the second factor.
Therefore, $|hk|_P|=|h|_P|$. By the same 
reason as $h$, $|hk|_P|=|hk|$.
By computation, 
$hk=\diag(a,a,1/a^2)$, which is hyperbolic.
We have just computed that $|\diag(a,a,1/a^2)|=
2\sqrt{6}\log|a|$.
To summarize, 
$$|g|=|h|=|h|_P|=|hk|_P|=|hk|=2\sqrt{6}\log|a|.$$
We finished the proof.
\qed


\begin{thebibliography}{99999}

\bibitem[B]{ballmann} 
W.~Ballmann,
\emph{Lectures on spaces of nonpositive curvature},
DMV-seminar, Band 25, Birkh\"{a}user, 1995.

\bibitem[BGS]{bagromsc}
W.~Ballmann, M.~Gromov, and V.~Schroeder,
\emph{Manifolds of nonpositive curvature},
Progress in Math., 61, Birkh\"{a}user, 1985.

\bibitem[BH]{bridhaef} 
M.~R.~Bridson and A.~Haefliger,
\emph{Metric spaces of non-positive curvature},
Grund. math. Wiss., Volume 319,
Springer-Verlag, 1999.

\bibitem[BBI]{burburiv}
D.~Burago, Yu.~ Burago, and S.~Ivanov,
\emph{A course in metric geometry},
Graduate Studies in Math., Volume 33, Amer. Math. Soc.,
2001.

\bibitem[Bu]{sebuyalo}
S.~V.~Buyalo,
\emph{Geodesics in Hadamard spaces},
St. Petersburg Math. J. {\bf 10} (1999), No. 2, 293--313.

\bibitem[E]{eberlein} 
P.~Eberlein,
\emph{Geometry of nonpositively curved manifolds},
Chicago Lectures in Math.,
The Univ. of Chicago Press, 1996.

\bibitem[F]{fedorchu}
V.~V.~Fedorchuk,
\emph{The fundamentals of dimension theory},
General Topology I,
Encycl. Math. Sci., Volume 17,
Springer-Verlag,
1990, pp. 91--192.

\bibitem[FSY]{fujshiya}
K.~Fujiwara, T.~Shioya, and S.~Yamagata,
\emph{Parabolic isometries of $\mathrm{CAT}(0)$ spaces
and $\mathrm{CAT}(0)$-dimension},
preprint, 2003.

\bibitem[J]{jurgjost}
J.~Jost,
\emph{Nonlinear Dirichlet forms},
New directions in Dirichlet forms,
AMS/IP Stud. Adv. Math., 8,
Amer. Math. Soc., 1998,
pp. 1--47.

\bibitem[K]{bkleiner}
B.~Kleiner,
\emph{The local structure of length spaces with curvature bounded above},
Math. Z. {\bf 231} (1999), 409--456.

\bibitem[LS1]{langsch1} 
U.~Lang and V.~Schroeder,
\emph{Kirszbraun's theorem and metric spaces of bounded curvature},
Geom. Funct. Anal. {\bf 7} (1997), 535--560.

\bibitem[LS2]{langsch2} 
\bysame,
\emph{Jung's theorem for Alexandrov spaces of curvature bounded above},
Ann. Global Anal. Geom. {\bf 15} (1997), 263--275.

\bibitem[L]{lytchak}
A.~Lytchak,
\emph{Rigidity of spherical buildings and joins},
preprint.

\bibitem[M]{uwemayer}
U.~Mayer,
\emph{Gradient flows on nonpositively curved metric spaces and
harmonic maps},
Comm. Anal. Geom., {\bf 6} (1998), Number 2, 199--253.

\end{thebibliography}
\end{document}